# Weyl Algebra Modules

Viktor Bekkert [*], Georgia Benkart[†], and Vyacheslav Futorny[‡]

**Abstract**

We investigate weight modules for finite and infinite Weyl algebras, classifying all such simple modules. We also study the representation type of the blocks of locally-finite weight module categories and describe indecomposable modules in tame blocks.

## 1  Introduction.

The $n$th *Weyl algebra* $A_n$ is the unital associative algebra over a field $\mathbb{K}$ with generators $x_i, \partial_i$, $i = 1, 2, \ldots, n$, which satisfy the defining relations

$$[x_i, x_j] = 0 = [\partial_i, \partial_j] \tag{1}$$
$$[\partial_i, x_j] = \delta_{i,j} 1, \tag{2}$$

for $1 \leq i, j \leq n$, expressed using the commutator $[a, b] = ab - ba$. Here $n$ may be a positive integer or countably infinite, and $\delta_{i,j}$ is the Kronecker delta.

Critical to our investigations of modules for $A_n$ is a realization of $A_n$ as a *generalized Weyl algebra* $D(\underline{\sigma}, \underline{a})$ of degree $n$ in the sense of [B1]. The ingredients for a generalized Weyl algebra $D(\underline{\sigma}, \underline{a})$ are a unital associative algebra $D$; an $n$-tuple of central elements $\underline{a} = (a_1, \ldots, a_n)$ of $D$; and an $n$-tuple of commuting automorphisms $\underline{\sigma} = (\sigma_1, \ldots, \sigma_n)$ of $D$ such that $\sigma_i(a_j) =$

[*]Supported in part by CNPq #301183/00-7.

[†]Supported in part by National Science Foundation Grants #DMS–9810361 (at MSRI) and #DMS–9970119.

[‡]Supported in part by CNPq #300679/1997-1, CRDF grant UM1-327, and National Science Foundation Grant #DMS–9810361 (at MSRI).

$a_j$ if $i \neq j$. The algebra $A = D(\underline{\sigma}, \underline{a})$ is generated over $D$ by elements $X_i$, $Y_i$, $i = 1, \ldots, n$, which satisfy the relations

$$X_i d = \sigma_i(d) X_i \qquad Y_i d = \sigma_i^{-1}(d) Y_i \qquad (3)$$
$$Y_i X_i = a_i \qquad X_i Y_i = \sigma_i(a_i) \qquad (4)$$
$$[X_i, X_j] = 0 = [Y_i, Y_j] \qquad 1 \leq i, j \leq n \qquad (5)$$
$$[Y_i, X_j] = 0 \qquad 1 \leq i \neq j \leq n, \qquad (6)$$

for all $d \in D$. (Here also $n$ may be a positive integer or countably infinite.)

The elements $t_i = \partial_i x_i$ in the Weyl algebra $A_n$ generate a polynomial algebra $D = \mathbb{K}[t_1, \ldots, t_n]$. By taking $\underline{a} = (a_1, \ldots, a_n)$ to be the tuple of elements $a_i = t_i$ in $D$, setting $X_i = x_i$, and $Y_i = \partial_i$, and letting $\underline{\sigma} = (\sigma_1, \ldots, \sigma_n)$ be the tuple of commuting automorphisms of $D$ given by $\sigma_i(t_j) = t_j - \delta_{i,j} 1$, where $\delta_{i,j} = 1$ if $i = j$ and is 0 otherwise, we obtain a realization of $A_n$ as the generalized Weyl algebra $D(\underline{\sigma}, \underline{a})$.

In [Bl], Block classified the simple modules for the Weyl algebra $A_1$ over an algebraically closed field $\mathbb{K}$ of characteristic 0 along with the simple modules for the Lie algebra $\mathfrak{sl}_2(\mathbb{K})$. An alternate approach using generalized Weyl algebras has been proposed in [B1], [B2], [BO1] when $\mathbb{K}$ is arbitrary. A module $V$ for a generalized Weyl algebra $D(\underline{\sigma}, \underline{a})$ is said to be a *weight module* if $V = \bigoplus_{\mathfrak{m} \in \mathfrak{max} D} V_\mathfrak{m}$, where $\mathfrak{max} D$ is the set of maximal ideals of $D$ and $V_\mathfrak{m} = \{v \in V \mid \mathfrak{m} v = 0\}$. Weight modules over $A_1$ (and over some other generalized Weyl algebras of degree 1) for arbitrary fields have been studied in [DGO]. For algebraically closed fields of characteristic 0, simple holonomic modules over the Weyl algebra $A_2$ have been classified in [BO2]. Weight and generalized weight modules for $A_n$ for $n < \infty$ (and for some other generalized Weyl algebras) over an algebraically closed field have been investigated in [BBL] and [BB]. In particular, [BBL] has given an explicit description of the weight modules for the complex Weyl algebras $A_n$ for $n < \infty$ and has used that to construct the weight modules having all weight spaces one-dimensional for the finite-dimensional simple complex Lie algebras. The paper [BB] has furnished a classification of the representation tame blocks in the category of locally-finite weight modules and described indecomposable modules in tame blocks. In this case all simple modules can be obtained as tensor products of simple modules over $A_1$.

In our paper we combine the techniques of [BB], [DGO], and [GP] to describe the simple weight $A_n$-modules for $n \leq \infty$ over *any* field. We also classify the blocks of tame representation type in the category of all locally-finite weight modules together with the indecomposable modules in each such block. From the classification of certain simple weight modules for $A_\infty$, we obtain many examples of simple $\mathbb{Z}$-graded modules for an infinite-dimensional Heisenberg Lie algebra with infinite-dimensional homogeneous components and a nonzero central charge.

Here is a brief outline of the paper. Section 2 introduces the main object of study – weight modules for Weyl algebras. In Section 3 we define the category $\mathcal{C}_\mathcal{O}$ and its skeleton. We determine all simple weight modules for Weyl algebras in Section 4. Various examples are described in Section 5. Using the classification of certain simple weight modules for $A_\infty$, we discuss simple $\mathbb{Z}$-graded modules with a nonzero central charge for the infinite-dimensional Heisenberg algebra in Section 6. Finally, in the last section we describe the blocks of tame representation type and classify indecomposable modules in those blocks.

## 2 Weight modules over Weyl algebras.

Assume for $1 \leq n \leq \infty$ that $A = A_n$, the $n$th Weyl algebra, and let $D = \mathbb{K}[t_1, \ldots, t_n]$ for $t_i = \partial_i x_i$ as above. Let $\mathcal{G}$ be the group generated by the automorphisms $\sigma_i$ ($1 \leq i \leq n$) of $D$, where $\sigma_i(t_j) = t_j - \delta_{i,j}1$. Then $\mathcal{G}$ acts on the set $\mathfrak{max}D$ of maximal ideals of $D$.

An $A$-module $V$ is a *weight module* if

$$V = \bigoplus_{\mathfrak{m} \in \mathfrak{max}D} V_\mathfrak{m}, \qquad V_\mathfrak{m} = \{v \in V \mid \mathfrak{m}v = 0\},$$

and if $V_\mathfrak{m} \neq 0$, we say that $\mathfrak{m}$ is a *weight* of $V$. The set $\mathfrak{supp}(V) = \{\mathfrak{m} \in \mathfrak{max}D \mid V_\mathfrak{m} \neq 0\}$ is the *support* of $V$. It is easy to see that $x_i V_\mathfrak{m} \subseteq V_{\sigma_i(\mathfrak{m})}$ and $\partial_i V_\mathfrak{m} \subseteq V_{\sigma_i^{-1}(\mathfrak{m})}$.

Let $\mathcal{W}(A)$ be the category of weight $A$-modules. For each subset $T \subset \mathfrak{max} D$ we will denote by $\mathcal{W}_T(A)$ the full subcategory in $\mathcal{W}(A)$ consisting of all modules $V$ with $\mathfrak{supp}(V) \subset T$. Each weight module $V$ can be decomposed into a direct sum of $A$-submodules:

$$V = \bigoplus_\mathcal{O} V_\mathcal{O}, \qquad V_\mathcal{O} := \bigoplus_{\mathfrak{m} \in \mathcal{O}} V_\mathfrak{m}$$

where $\mathcal{O}$ runs over the orbits of $\mathcal{G}$ on $\mathfrak{max}D$. Hence, the category $\mathcal{W}(A)$ decomposes into a sum of full subcategories corresponding to the orbits of $\mathcal{G}$. In particular, if $V$ is indecomposable, then its support belongs to a single orbit $\mathcal{O}$.

An orbit is said to be *cyclic* (resp. *linear*) if $\tau(\mathfrak{m}) = \mathfrak{m}$ for some $\tau \in \mathcal{G}$, $\tau \neq 1$ (resp. $\tau(\mathfrak{m}) \neq \mathfrak{m}$ for all $\tau \in \mathcal{G}$, $\tau \neq 1$). It is evident that if $\mathrm{char}\,\mathbb{K} = 0$, then each orbit is linear. On the contrary, if $\mathrm{char}\,\mathbb{K} = p > 0$, then each orbit $\mathcal{O}$ is cyclic and $|\mathcal{O}| < \infty$ if $n < \infty$.

**Remark 2.1.** *Let $V \in \mathcal{W}(A)$ be simple and assume $\mathfrak{m} \in \mathfrak{supp}(V)$. Suppose that the corresponding orbit is linear. Following [BO1], we consider $A(\mathfrak{m}) = A/A\mathfrak{m} \in \mathcal{W}(A)$ and the maximal submodule $N(\mathfrak{m}) \subset A(\mathfrak{m})$ which is the sum of all submodules that intersect $D/\mathfrak{m}$ trivially. Then $V \cong A(\mathfrak{m})/N(\mathfrak{m})$. Hence, for linear orbits the problem is to describe the submodule $N(\mathfrak{m})$.*

Assume now that $\mathrm{char}\,\mathbb{K} = p > 0$. The algebra $A = A_n$ is $\mathbb{Z}^n$-graded: $A = \sum_{\underline{j} \in \mathbb{Z}^n} A_{\underline{j}}$. For $\underline{j} = (j_1, \ldots, j_n) \in \mathbb{Z}^n$, the homogeneous space $A_{\underline{j}} = Dv_{\underline{j}}$, where $v_{\underline{j}} = v_{j_1}(1) \ldots v_{j_n}(n)$, and $v_j(i) = x_i^j$, $v_{-j}(i) = \partial_i^j$, $j > 0$, $v_0(i) = 1$ for all $i = 1, \ldots, n$. We define the *Veronese subalgebra* $A_{[p]} = \sum_{\underline{j} \in \mathbb{Z}^n} A_{p\underline{j}}$ (cf. [BO1]). If $\mathfrak{m}$ is a maximal ideal of $D$, then $A_{[p]}\mathfrak{m} \subseteq A_{[p]}$ as $D \subseteq A_{\underline{0}}$.

**Remark 2.2.** *Let $\mathrm{char}\,\mathbb{K} = p > 0$ and $V \in \mathcal{W}(A)$ be a simple module for $A = A_n$. Assume $V_\mathfrak{m} \neq 0$ for some maximal ideal $\mathfrak{m}$ of $D$. Define an $A$-module by $A(\mathfrak{m}) = A \otimes_{A_{[p]}} V_\mathfrak{m}$. Then $V \cong A(\mathfrak{m})/N(\mathfrak{m})$ where $N(\mathfrak{m})$ is the maximal submodule in $A(\mathfrak{m})$ which trivially intersects $D/\mathfrak{m}$. The problem of describing such simple modules $V$ reduces to determining the simple modules for $A_{[p]}/A_{[p]}\mathfrak{m}$ and the maximal submodule $N(\mathfrak{m})$.*

Remarks 2.1 and 2.2 provide the construction of all the simple weight $A_n$-modules, reducing the problem of their classification to the description of the maximal submodule $N(\mathfrak{m})$. Because determining the maximal submodule $N(\mathfrak{m})$ can be highly nontrivial, in what follows we will adopt an alternate approach based on the ideas of [DGO] and [BB].

**2.3.** We say that the maximal ideal $\mathfrak{m}$ of $D$ is a *break with respect to $i$* if $t_i \in \mathfrak{m}$ for some $i \in \{1, \ldots, n\}$. An orbit $\mathcal{O}$ is *degenerate with respect to $i$* if it contains a break with respect to $i$ (for some $\mathfrak{m} \in \mathcal{O}$). Often we simply say that $\mathcal{O}$ is degenerate without specifying $i$ or $\mathfrak{m}$.

A maximal ideal $\mathfrak{m}$ of $D$ is a *maximal break with respect to $\underline{I} = \underline{I}(\mathfrak{m}) \subset \{1, \ldots, n\}$* if $t_i \in \mathfrak{m}$ for each $i \in \underline{I}$, and $t_j \notin \tau(\mathfrak{m})$ for each $j \notin \underline{I}$ and each $\tau \in \mathcal{G}$. We set $\underline{I}^c := \{1, \ldots, n\} \setminus \underline{I}$ and say that the maximal break has order $|\underline{I}|$ (which may be infinite if $n = \infty$).

To see that maximal breaks exist for degenerate orbits when $n < \infty$, suppose $\mathcal{O} = \mathcal{O}(\mathfrak{m})$ is degenerate. If $\mathfrak{n} \in \mathcal{O}$ is a break with respect to $i$, then $t_i \in \mathfrak{n}$. Moreover, for any product $\tau = \sigma_1^{k_1} \cdots \sigma_n^{k_n} \in \mathcal{G}$ with the $\sigma_i$ term omitted, we have $t_i \in \tau(\mathfrak{n})$ as $\sigma_j(t_i) = t_i$ for $j \neq i$. Suppose $\underline{I} = \{i_1, \ldots, i_s\} \subseteq \{1, \ldots, n\}$ is a maximal set of indices such that $\mathcal{O}$ has breaks $\mathfrak{n}_{i_1}, \ldots, \mathfrak{n}_{i_s}$ relative to $i_1, \ldots, i_s$ respectively and no breaks with respect to any $j \in \underline{I}^c$. Assume $\mathfrak{n}_i = \sigma_1^{r_{1,i}} \cdots \sigma_n^{r_{n,i}}(\mathfrak{m})$ for each $i \in \underline{I}$. Then for $i \in \underline{I}$, $t_i \in \sigma_i^{r_{i,i}}(\mathfrak{m})$, hence $t_i \in \mathfrak{n} := \sigma_{i_1}^{r_{i_1,i_1}} \cdots \sigma_{i_s}^{r_{i_s,i_s}}(\mathfrak{m})$, and $\mathfrak{n}$ is a maximal break in $\mathcal{O}$ with respect to $\underline{I} = \{i_1, \ldots, i_s\}$.

**2.4.** Because it is not apparent that maximal breaks exist for $A_\infty$, we will assume that the degenerate orbits $\mathcal{O}$ we consider for $A_\infty$ have a maximal break $\mathfrak{m}$. In this case, the break set $\underline{I}(\mathfrak{m})$ may be finite or infinite.

**2.5.** In what follows, we let $\underline{\tau}(\mathfrak{m}) = (\tau_1, \cdots, \tau_n)$, where $\tau_i = \sigma_i$ if $\sigma_i(\mathfrak{m}) = \mathfrak{m}$ and $\tau_i = 1$ otherwise. Each such $\tau_i$ induces an automorphism on $D/\mathfrak{m}$.

## 3 Category $\mathcal{C}_\mathcal{O}$ and its skeleton.

We begin this section with the abstract concepts needed to describe Weyl algebra modules, which can be found in ([GR], Chap. 2). Ultimately we specialize to very particular categories determined from these modules.

### 3.1. Categories

Let $\mathbb{F}$ be a field. A category $\mathcal{C}$ is said to be an $\mathbb{F}$-category if each morphism set $\mathcal{C}(\alpha, \beta)$ is equipped with an $\mathbb{F}$-bimodule structure, the composition is $\mathbb{F}$-linear with respect to both the left and right $\mathbb{F}$-module structures, and $(\alpha\lambda)\beta = \alpha(\lambda\beta)$ for any possible morphisms $\alpha, \beta$ and $\lambda \in \mathbb{F}$. If additionally $\alpha\lambda = \lambda\alpha$ for all morphisms $\alpha$ and $\lambda \in \mathbb{F}$, then we say that $\mathcal{C}$ is an $\mathbb{F}$-linear category. An $\mathbb{F}$-algebra $A$ gives rise to an $\mathbb{F}$-category with just one object - say $\alpha$ such that $\mathcal{C}(\alpha, \alpha) = A$, and the composition is multiplication in $A$.

With any category $\mathcal{C}$, we can associate an $\mathbb{F}$-linear category $\mathbb{F}\mathcal{C}$ whose objects are the same $\mathfrak{ob}(\mathcal{C}) = \mathfrak{ob}(\mathbb{F}\mathcal{C})$, and the morphism space $\mathbb{F}\mathcal{C}(\alpha, \beta)$ has as its basis over $\mathbb{F}$ the elements of $\mathcal{C}(\alpha, \beta)$. The composition in $\mathbb{F}\mathcal{C}$ is the $\mathbb{F}$-linear extension of the composition in $\mathcal{C}$.

An $\mathbb{F}$-linear functor between two $\mathbb{F}$-linear categories $\mathcal{C}$ and $\mathcal{D}$ is a functor $F : \mathcal{C} \to \mathcal{D}$ whose defining maps $F(\alpha, \beta) : \mathcal{C}(\alpha, \beta) \to \mathcal{D}(F\alpha, F\beta)$ are $\mathbb{F}$-linear for all $\alpha, \beta \in \mathfrak{ob}\mathcal{C}$.

Let $\mathbb{F} - \mathfrak{mod}$ denote the category of $\mathbb{F}$-vector spaces and let $\mathfrak{Ab}$ be the category of abelian groups. Given an $\mathbb{F}$-category $\mathcal{C}$ we denote by $\mathcal{C}\text{-}\mathfrak{mod}$ the category of all additive functors $M : \mathcal{C} \to \mathfrak{Ab}$. The functors $M$ are called $\mathcal{C}$-*modules*, or more precisely *left $\mathcal{C}$-modules*. For an $\alpha \in \mathfrak{ob}\mathcal{C}$, the elements of $M(\alpha)$ are the *elements of the module $M$* (at $\alpha$). To make the action more module-like, we shall write $au$ instead of $M(a)u$ for $u \in M(\alpha)$ and $a \in \mathcal{C}(\alpha, \beta)$. For each $\alpha \in \mathfrak{ob}\mathcal{C}$, the group $M(\alpha)$ becomes an $\mathbb{F}$-vector space if we put $\lambda u = (\lambda 1_\alpha)u$ for all $u \in M(\alpha)$ and $\lambda \in \mathbb{F}$. By $\mathcal{C}$-$\mathfrak{fdmod}$ we mean the full subcategory of all locally finite-dimensional objects $M$ in $\mathcal{C}$-$\mathfrak{mod}$, that is $\dim_\mathbb{F} M(\alpha) < \infty$ for all $\alpha \in \mathfrak{ob}\mathcal{C}$.

If $N$ is a $\mathcal{C}$-*submodule* of $M$, then $N(\alpha)$ is a subspace of $M(\alpha)$ for all $\alpha \in \mathfrak{ob}\mathcal{C}$ and $au \in N(\beta)$ for all $a \in \mathcal{C}(\alpha, \beta)$ and $u \in N(\alpha)$. The module $M$ is *simple* if it has no nontrivial $\mathcal{C}$-submodules; while $M = M_1 \oplus M_2$ if $M_1$ and $M_2$ are submodules such that $M(\alpha) = M_1(\alpha) \oplus M_2(\alpha)$ for all $\alpha \in \mathfrak{ob}\mathcal{C}$. The module $M$ is *indecomposable* if it cannot be written as a direct sum $M = M_1 \oplus M_2$ of proper submodules.

Right $\mathcal{C}$-modules may be defined as $\mathcal{C}^{\mathfrak{op}}$-modules, where $\mathcal{C}^{\mathfrak{op}}$ is the category opposite to $\mathcal{C}$. In this case we write $va$ for $v \in M(\beta)$ and for $a \in \mathcal{C}(\alpha, \beta) = \mathcal{C}^{\mathfrak{op}}(\beta, \alpha)$. The category of right $\mathcal{C}$-modules will be denoted by $\mathfrak{mod}$-$\mathcal{C}$ and the full subcategory of locally finite-dimensional objects by $\mathfrak{fdmod}$-$\mathcal{C}$.

For $\mathbb{F}$-categories $\mathcal{C}$ and $\mathcal{D}$ a $\mathcal{C} - \mathcal{D}$-*bimodule* is an additive functor $B : \mathcal{C} \times \mathcal{D}^{\mathfrak{op}} \to \mathfrak{Ab}$. If $u \in B(\gamma, \delta)$, then $u$ is an element of the module $B$ with *source* $\gamma$ and *target* $\delta$, and we write $aub$ rather than $B(a,b)u$ for $a \in \mathcal{C}(\gamma, \gamma')$ and $b \in \mathcal{D}(\delta', \delta)$. Any $\mathbb{F}$-category $\mathcal{C}$ can be viewed as a $\mathcal{C} - \mathcal{C}$-bimodule mapping a pair of objects $(\beta, \gamma)$ to $\mathcal{C}(\beta, \gamma)$.

A category $\mathcal{C}$ is said to be *basic* if
- all its objects are pairwise nonisomorphic;
- for each object $\alpha$ there are no nontrivial idempotents in $\mathcal{C}(\alpha, \alpha)$.

A full subcategory $\mathcal{S}$ is a *skeleton* of a category $\mathcal{C}$ if it is basic, and each object $\alpha \in \mathfrak{ob}\mathcal{C}$ is isomorphic to a direct summand of a (finite) direct sum of some objects of $\mathcal{S}$. If the category $\mathcal{C}$ has a unique direct decomposition property, then it has a skeleton which is unique up to isomorphism. The natural inclusion functor $\mathcal{I}: \mathcal{S} \to \mathcal{C}$ of a skeleton $\mathcal{S}$ into $\mathcal{C}$, is an equivalence of categories. By this functor, $\mathcal{C}$ becomes a $\mathcal{C} - \mathcal{S}$-bimodule in an obvious way. Tensoring this bimodule over $\mathcal{S}$ furnishes equivalences $\mathcal{S}$-$\mathfrak{mod} \to \mathcal{C}$-$\mathfrak{mod}$ and $\mathcal{S}$-$\mathfrak{fdmod} \to \mathcal{C}$-$\mathfrak{fdmod}$. This is a reformulation of Morita equivalence in the categorical context.

### 3.2. Quivers

A *quiver* $Q$ is a tuple $(Q_0, Q_1, \mathfrak{s}, \mathfrak{e})$ consisting of a set $Q_0$ of *vertices*, a set $Q_1$ of arrows, and maps $\mathfrak{s}, \mathfrak{e} : Q_1 \to Q_0$ which specify the starting and ending vertices. A *path* $p$ in $Q$ of length $\ell(p) = n \geq 1$ is a sequence of arrows $a_n, \ldots, a_1$ such that $\mathfrak{s}(a_{i+1}) = \mathfrak{e}(a_i)$ for $1 \leq i < n$. Set $\mathfrak{s}(p) = \mathfrak{s}(a_1)$ and $\mathfrak{e}(p) = \mathfrak{e}(a_n)$. Then the concatenation $p'p$ of two paths $p$, $p'$ is defined in the natural way whenever $\mathfrak{s}(p') = \mathfrak{e}(p)$. Every vertex $a \in Q_0$ determines a path $1_a$ (of length 0) with $\mathfrak{s}(1_a) = a$ and $\mathfrak{e}(1_a) = a$. A quiver $Q$ determines a category $\mathfrak{path}Q$ with objects the vertices of $Q$ and morphisms from vertex $a$ to vertex $b$ being the paths from $a$ to $b$. The composition in $\mathfrak{path}Q$ of paths of positive length is just concatenation, and the path $1_a$ acts as the identity on all paths for which composition makes sense. The corresponding path algebra, which we will abbreviate simply $\mathbb{F}Q$, has an $\mathbb{F}$-basis consisting of the paths of $Q$ with multiplication given by concatenation of paths. If $\mathcal{R}$ is the ideal of $\mathbb{F}Q$ generated by a family $\{\rho_i\}$ of morphisms, we say that $\mathbb{F}Q/\mathcal{R}$ is the $\mathbb{F}$-*linear category defined by the quiver $Q$ and the relations $\{\rho_i\}$.*

### 3.3. Category $\mathcal{C}_\mathcal{O}$

In this section we prove that the category $\mathcal{W}_\mathcal{O}(A)$ of weight modules for the Weyl algebra $A = A_n$ ($1 \leq n \leq \infty$) having support in the orbit $\mathcal{O}$, (which is assumed to have a maximal break if $\mathcal{O}$ is degenerate and $n = \infty$), is equivalent to $\mathcal{C}_\mathcal{O}$-$\mathfrak{mod}$, the category of modules over a certain category $\mathcal{C}_\mathcal{O}$ first introduced in [DGO].

Assume $\mathfrak{m}$ is a fixed maximal break of the orbit $\mathcal{O}$ if $\mathcal{O}$ is degenerate and is any fixed element of $\mathcal{O}$ otherwise. For a given $\mathfrak{n} \in \mathcal{O}$ we denote by $\sigma_\mathfrak{n}$ the element of $\mathcal{G}/\mathfrak{stab}(\mathfrak{m})$ (where $\mathfrak{stab}(\mathfrak{m})$ is the stabilizer subgroup of $\mathfrak{m}$ in $\mathcal{G}$) such that $\mathfrak{n} = \sigma_\mathfrak{n}(\mathfrak{m})$. Thus, $\sigma_\mathfrak{n}$ induces an isomorphism $D/\mathfrak{m} \to D/\mathfrak{n}$, which we again denote $\sigma_\mathfrak{n}$, and its inverse induces an isomorphism $\sigma_\mathfrak{n}^{-1} : D/\mathfrak{n} \to D/\mathfrak{m}$.

We define $\mathcal{C}_\mathcal{O}$ as the $D/\mathfrak{m}$-category with $\mathfrak{ob}\mathcal{C}_\mathcal{O} = \mathcal{O}$, generated over $D/\mathfrak{m}$ by the set of morphisms $\{X_{\mathfrak{n},i}, Y_{\mathfrak{n},i} \mid \mathfrak{n} \in \mathcal{O}, 1 \leq i \leq n\}$, where $X_{\mathfrak{n},i} : \mathfrak{n} \to \sigma_i(\mathfrak{n})$ and $Y_{\mathfrak{n},i} : \sigma_i(\mathfrak{n}) \to \mathfrak{n}$, subject to the relations:

- $X_{\mathfrak{n},i}\lambda = \lambda X_{\mathfrak{n},i}$, $Y_{\mathfrak{n},i}\lambda = \lambda Y_{\mathfrak{n},i}$, $Y_{\mathfrak{n},i}X_{\mathfrak{n},i} = \sigma_\mathfrak{n}^{-1}(\bar{t}_i)1_\mathfrak{n}$ and $X_{\mathfrak{n},i}Y_{\mathfrak{n},i} = \sigma_\mathfrak{n}^{-1}(\bar{t}_i)1_{\sigma_i(\mathfrak{n})}$, $\bar{t}_i = t_i + \mathfrak{n}$, for each $\lambda \in D/\mathfrak{m}$, $\mathfrak{n} \in \mathcal{O}$, and each $1 \leq i \leq n$ such that $\sigma_i(\mathfrak{m}) \neq \mathfrak{m}$;

- $X_{\mathfrak{n},i}\lambda = \sigma_i(\lambda)X_{\mathfrak{n},i}$, $Y_{\mathfrak{n},i}\sigma_i(\lambda) = \lambda Y_{\mathfrak{n},i}$, $Y_{\mathfrak{n},i}X_{\mathfrak{n},i} = \sigma_\mathfrak{n}^{-1}(\bar{t}_i)1_\mathfrak{n}$ and $X_{\mathfrak{n},i}Y_{\mathfrak{n},i} = \sigma_i\sigma_\mathfrak{n}^{-1}(\bar{t}_i)1_{\sigma_i(\mathfrak{n})}$, $\bar{t}_i = t_i + \mathfrak{n}$, for each $\lambda \in D/\mathfrak{m}$, $\mathfrak{n} \in \mathcal{O}$, and each $1 \leq i \leq n$ such that $\sigma_i(\mathfrak{m}) = \mathfrak{m}$;

- $U_{\mathfrak{n},i}V_{\mathfrak{p},j} - V_{\mathfrak{q},j}U_{\mathfrak{r},i} = 0$ for all $j \neq i$ and all possible $U, V \in \{X, Y\}$, $\mathfrak{n}, \mathfrak{p}, \mathfrak{q}, \mathfrak{r} \in \mathcal{O}$, for which the last equality makes sense.

Note that the category $\mathcal{C}_\mathcal{O}$ is $\mathbb{F}$-linear when $\sigma_i(\mathfrak{m}) \neq \mathfrak{m}$ for all $i$; in particular, $\mathcal{C}_\mathcal{O}$ is always $\mathbb{F}$-linear when $\mathbb{F}$ has characteristic 0.

**Proposition 3.4.** *Let $A = A_n$, the nth Weyl algebra over the field $\mathbb{K}$, and let $\mathcal{O}$ be an orbit of $\mathfrak{max}D$ (which is assumed to have a maximal break if $\mathcal{O}$ is degenerate and $n = \infty$). Then $\mathcal{W}_{\mathcal{O}}(A) \cong \mathcal{C}_{\mathcal{O}}\text{-}\mathfrak{mod}$.*

**Proof.** The proof is similar to that of Proposition 2.2 in [DGO]. We assume $\mathfrak{m}$ is the designated maximal ideal of $\mathcal{O}$. Let $V = \bigoplus_{\mathfrak{n} \in \mathcal{O}} V_{\mathfrak{n}}$ belong to $\mathcal{W}_{\mathcal{O}}(A)$. For each $\mathfrak{n} \in \mathcal{O}$ set $M_V(\mathfrak{n}) = V_{\mathfrak{n}}$. Using the isomorphism $\sigma_{\mathfrak{n}} : D/\mathfrak{m} \to D/\mathfrak{n}$, we can view $M_V(\mathfrak{n})$ as a $D/\mathfrak{m}$-vector space via $\bar{d}v := \sigma_{\mathfrak{n}}(\bar{d})v$, $(\bar{d} = d + \mathfrak{m})$. For $v \in M_V(\mathfrak{n})$ and $w \in M_V(\sigma_i(\mathfrak{n}))$, we define $X_{\mathfrak{n},i}v := x_i v \in M_V(\sigma_i(\mathfrak{n}))$ and $Y_{\mathfrak{n},i}w := \partial_i w \in M_V(\mathfrak{n})$. Then if $\bar{d} \in D/\mathfrak{m}$, we have

$$\begin{aligned} X_{\mathfrak{n},i}\bar{d}v &= x_i \sigma_{\mathfrak{n}}(\bar{d})v = \sigma_i(\sigma_{\mathfrak{n}}(\bar{d}))x_i v \\ &= \begin{cases} \bar{d}X_{\mathfrak{n},i}v & \text{if } \sigma_i(\mathfrak{m}) \neq \mathfrak{m} \\ \sigma_i(\bar{d})X_{\mathfrak{n},i}v & \text{if } \sigma_i(\mathfrak{m}) = \mathfrak{m}, \end{cases} \end{aligned}$$

as $\sigma_i(\mathfrak{n}) = \mathfrak{n}$ for all $\mathfrak{n} \in \mathcal{O}$, and $\sigma_i(\sigma_{\mathfrak{n}}(\bar{d})) = \sigma_{\mathfrak{n}}(\sigma_i(\bar{d}))$, whenever $\sigma_i(\mathfrak{m}) = \mathfrak{m}$. Likewise $Y_{\mathfrak{n},i}\bar{d}w = \bar{d}Y_{\mathfrak{n},i}w$. We also have $Y_{\mathfrak{n},i}X_{\mathfrak{n},i}v = \partial_i x_i v = t_i v = (t_i + \mathfrak{n})v = \sigma_{\mathfrak{n}}^{-1}(\bar{t}_i)v$ for $\bar{t}_i = t_i + \mathfrak{n}$, and $X_{\mathfrak{n},i}Y_{\mathfrak{n},i}w = x_i \partial_i w = \sigma_i(t_i)w = (\sigma_i(t_i) + \sigma_i(\mathfrak{n}))w = \sigma_{\mathfrak{n}}^{-1}(\bar{t}_i)w$. The reason for the last equality is that there is an action of $\mathcal{G}$ on $\mathcal{G}/\mathfrak{stab}(\mathfrak{m})$, and under this action $\sigma_i \sigma_{\mathfrak{n}} = \sigma_{\sigma_i(\mathfrak{n})}$. Hence, $M_V$ ($M_V : \mathcal{C}_{\mathcal{O}} \to D/\mathfrak{m}\text{-}\mathfrak{mod}$) is a $\mathcal{C}_{\mathcal{O}}$-module, and we have a functor

$$F : \mathcal{W}_{\mathcal{O}}(A) \to \mathcal{C}_{\mathcal{O}}\text{-}\mathfrak{mod}, \qquad V \mapsto M_V. \tag{7}$$

Conversely, for each $M \in \mathcal{C}_{\mathcal{O}}\text{-}\mathfrak{mod}$, we define the $A$-module $V_M := \oplus_{\mathfrak{n} \in \mathcal{O}} M(\mathfrak{n})$ where $dv := \sigma_{\mathfrak{n}}^{-1}(\bar{d})v$, $\bar{d} = d + \mathfrak{n} \in D/\mathfrak{n}$, $x_i v := X_{\mathfrak{n},i}v$, and $\partial_i v := Y_{\sigma_i^{-1}(\mathfrak{n}),i}v$ for $v \in M(\mathfrak{n})$. This gives a functor

$$F' : \mathcal{C}_{\mathcal{O}}\text{-}\mathfrak{mod} \to \mathcal{W}_{\mathcal{O}}(A), \qquad M \mapsto V_M \tag{8}$$

which is inverse to $F$. $\square$

### 3.5. Skeleton of $\mathcal{C}_{\mathcal{O}}$

For a given orbit $\mathcal{O}$ we define the set $\mathfrak{B}_{\mathcal{O}}$ according to the following rule: If $\mathcal{O}$ is nondegenerate, then $\mathfrak{B}_{\mathcal{O}} = \{\mathfrak{m}\}$ for the fixed $\mathfrak{m} \in \mathcal{O}$ used to define $\mathcal{C}_{\mathcal{O}}$, and $\underline{I}(\mathfrak{m}) = \emptyset$. If $\mathcal{O}$ is a linear degenerate orbit, then $\mathfrak{m} \in \mathcal{O}$ is the fixed maximal break used in the construction of $\mathcal{C}_\mathcal{O}$. We assume that the set of breaks is $\underline{I} = \underline{I}(\mathfrak{m}) = \{i_1, \ldots, i_s\}$. We set $\mathfrak{B}_\mathcal{O} = \{\sigma_{i_1}^{\delta_1} \cdots \sigma_{i_s}^{\delta_s}(\mathfrak{m}) \,|\, \delta_j \in \{0,1\}\}$ (where only finitely many $\delta_j$ are nonzero when $n = \infty$). In this case, for each maximal ideal $\mathfrak{n} = \sigma_{i_1}^{\delta_1} \cdots \sigma_{i_s}^{\delta_s}(\mathfrak{m}) \in \mathfrak{B}_\mathcal{O}$ we set

$$\mathcal{O}_\mathfrak{n} := \{\sigma_1^{\gamma_1} \cdots \sigma_n^{\gamma_n}(\mathfrak{n}) \,|\, \gamma_i = (-1)^{\delta_i+1}k,\ k \in \mathbb{N} \text{ if } i \in \underline{I},\text{ and } \gamma_i \in \mathbb{Z} \text{ if } i \in \underline{I}^c\}, \tag{9}$$

(where only finitely many $\gamma_i$ are nonzero when $n = \infty$). If $\mathcal{O}$ is a cyclic degenerate orbit, and $\mathfrak{m} \in \mathcal{O}$ is the maximal break used to define $\mathcal{C}_\mathcal{O}$, then $\mathfrak{B}_\mathcal{O} := \{\mathfrak{m}\}$ and $\mathcal{O}_\mathfrak{m} := \mathcal{O}$.

**3.6.** We introduce an equivalence relation $\sim$ on the set of maximal ideals $\mathfrak{n} \in \mathfrak{max}D$. This relation is the transitive extension of the relation specified by $\mathfrak{n} \sim \sigma_i(\mathfrak{n})$ if and only if $t_i \notin \mathfrak{n}$.

**Lemma 3.7.** *Assume $\mathfrak{n}$ and $\mathfrak{p}$ belong to $\mathcal{O}$. Then $\mathfrak{n} \sim \mathfrak{p}$ if and only if $\mathfrak{n}$ and $\mathfrak{p}$ are isomorphic in $\mathcal{C}_\mathcal{O}$.*

**Proof.** Let $\mathfrak{m}$ be the fixed maximal break of $\mathcal{O}$ if $\mathcal{O}$ is degenerate and be the unique element of $\mathfrak{B}_\mathcal{O}$ if $\mathcal{O}$ is nondegenerate. We assume for each $\mathfrak{n} \in \mathcal{O}$ that $\sigma_\mathfrak{n} \in \mathcal{G}/\mathfrak{stab}(\mathfrak{m})$ satisfies $\sigma_\mathfrak{n}(\mathfrak{m}) = \mathfrak{n}$ as in the proof of Proposition 3.4.

It suffices to argue for any $\mathfrak{n} \in \mathcal{O}$ that $\mathfrak{n} \sim \sigma_i(\mathfrak{n})$ if and only if $\mathfrak{n}$ and $\sigma_i(\mathfrak{n})$ are isomorphic in $\mathcal{C}_\mathcal{O}$.

($\Longrightarrow$) Assume $\mathfrak{n} \sim \sigma_i(\mathfrak{n})$. Then $t_i \notin \mathfrak{n}$ so $t_i + \mathfrak{n} \neq 0$. As $\sigma_\mathfrak{n}^{-1} : D/\mathfrak{n} \to D/\mathfrak{m}$ is an isomorphism, $\sigma_\mathfrak{n}^{-1}(t_i + \mathfrak{n}) \neq 0$ and by (3.3), $X_{\mathfrak{n},i}$ and $Y_{\mathfrak{n},i}$ are invertible in $\mathcal{C}_\mathcal{O}$. Hence $\mathfrak{n}$ and $\sigma_i(\mathfrak{n})$ are isomorphic in $\mathcal{C}_\mathcal{O}$.

($\Longleftarrow$) Assume $\mathfrak{n}$ and $\sigma_i(\mathfrak{n})$ are isomorphic in $\mathcal{C}_\mathcal{O}$. It follows from the definition of the morphisms in $\mathcal{C}_\mathcal{O}$ that $\mathcal{C}_\mathcal{O}(\mathfrak{n}, \sigma_i(\mathfrak{n})) = (D/\mathfrak{m})X_{\mathfrak{n},i}$ and $\mathcal{C}_\mathcal{O}(\sigma_i(\mathfrak{n}), \mathfrak{n}) = (D/\mathfrak{m})Y_{\mathfrak{n},i}$ for $1 \leq i \leq n$. Hence $X_{\mathfrak{n},i}$ and $Y_{\mathfrak{n},i}$ are isomorphisms. If $t_i \in \mathfrak{n}$, then $\sigma_\mathfrak{n}^{-1}(t_i + \mathfrak{n}) = 0 \in D/\mathfrak{m}$ and $Y_{\mathfrak{n},i}X_{\mathfrak{n},i} = \sigma_\mathfrak{n}^{-1}(t_i + \mathfrak{n})1_\mathfrak{n} = 0$, a contradiction. □

**Remark 3.8.** For $\mathfrak{n} \in \mathfrak{B}_\mathcal{O}$, the equivalence class of $\mathfrak{n}$ is exactly $\mathcal{O}_\mathfrak{n}$, and $\mathcal{O}$ is the disjoint union of the sets $\mathcal{O}_\mathfrak{n}$, $\mathfrak{n} \in \mathfrak{B}_\mathcal{O}$.

**Corollary 3.9.** *The full subcategory $\mathcal{S}_\mathcal{O}$ with $\mathfrak{ob}\mathcal{S}_\mathcal{O} = \mathfrak{B}_\mathcal{O}$ is a skeleton of the category $\mathcal{C}_\mathcal{O}$.*

**Proof.** Assume $\mathfrak{m}$ is the designated maximal ideal in $\mathcal{O}$, and the corresponding set of breaks is $\underline{I} = \underline{I}(\mathfrak{m}) = \{i_1, \ldots, i_s\}$ (which is empty if $\mathcal{O}$ is nondegenerate). In the degenerate case, the maximal ideal $\sigma_{i_1}^{\delta_1} \cdots \sigma_{i_s}^{\delta_s}(\mathfrak{m})$ has breaks at $\{i_k \mid \delta_k = 0\}$. Thus, these ideals are pairwise nonisomorphic. Moreover, if $\mathfrak{n} = \sigma_1^{r_1} \cdots \sigma_n^{r_n}(\mathfrak{m}) \in \mathcal{O}$ (with only finitely many $r_i$ nonzero when $n = \infty$), then by Lemma 3.7, $\mathfrak{n}$ is isomorphic to $\sigma_{i_1}^{\delta_1} \cdots \sigma_{i_s}^{\delta_s}(\mathfrak{m}) \in \mathfrak{B}_\mathcal{O}$ such that $\delta_i = 1$ whenever $i \in \underline{I}$ and $r_i \geq 1$, and $\delta_i = 0$ otherwise. In the nondegenerate case, $\mathfrak{B}_\mathcal{O} = \{\mathfrak{m}\}$ and every element of $\mathcal{O}$ is isomorphic to $\mathfrak{m}$ by the lemma. Consequently, in both cases $\mathfrak{B}_\mathcal{O}$ is a skeleton of $\mathcal{C}_\mathcal{O}$. □

### 3.10. Algebraic description of the skeleton

For a given field $\mathbb{F}$ and an arbitrary subset $I$ of the set $\mathbb{N}$ of positive integers, we define the category $\mathcal{A} = \mathcal{A}(\mathbb{F}, I)$ as the $\mathbb{F}$-linear category with the set of objects $\mathfrak{ob}\mathcal{A} := \{0, 1\}^{|I|}$ (which we assume have only finitely many nonzero components when $|I| = \infty$) generated (over $\mathbb{F}$) by the set of morphisms, $\mathcal{A}_1 := \{a_{\alpha,i}, b_{\alpha,i} \mid \alpha \in \mathfrak{ob}\mathcal{A}, \ i \in I, \text{ and } \alpha_i = 0\}$, where $a_{\alpha,i} : \alpha \to \beta$ and $b_{\alpha,i} : \beta \to \alpha$, such that $\beta_j = \alpha_j$ for all $j \neq i$ and $\beta_i = 1$, subject to the relations:

- $a_{\alpha,i} b_{\alpha,i} = b_{\alpha,i} a_{\alpha,i} = 0$ for each $a_{\alpha,i}, b_{\alpha,i} \in \mathcal{A}_1$;

- $u_{\alpha,i} v_{\beta,j} - v_{\gamma,j} u_{\delta,i} = 0$ for all $j \neq i$ and all possible $u, v \in \{a, b\}$, $\alpha, \beta, \gamma, \delta \in \mathfrak{ob}\mathcal{A}$, for which the last equality makes sense.

When $I$ is empty, let $\mathcal{A}(\mathbb{F}, \emptyset)$ be the category with a unique object, say $\omega$, and with morphism set $\mathbb{F}1_\omega$.

The $\mathbb{F}$-algebra corresponding to the category $\mathcal{A}(\mathbb{F}, I)$ above consists of finite $\mathbb{F}$-linear combinations of morphisms in the category, and the product is simply composition of morphisms whenever it is defined and is 0 otherwise. It has a unit element if $I$ is finite. We adopt the same notation $\mathcal{A}(\mathbb{F}, I)$ for the algebra, as it will be evident from the context which one is meant. The algebra $\mathcal{A}(\mathbb{F}, I)$ is finite-dimensional when $|I| = s < \infty$, since for $I = \{i_1, \cdots, i_s\}$, $\mathcal{A}(\mathbb{F}, I) \cong \mathcal{A}(\mathbb{F}, \{i_1\}) \otimes \cdots \otimes \mathcal{A}(\mathbb{F}, \{i_s\})$. It is easy to see that algebra $\mathcal{A}(\mathbb{F}, \{i\})$ is isomorphic to the algebra $\overline{Q}_1 := \mathbb{F}\mathcal{Q}_1/\mathcal{R}$ corresponding to the following quiver and relations:

$$\mathcal{Q}_1 : \quad \underset{b \quad 2}{\overset{1 \quad a}{\circ \rightleftarrows \circ}} \qquad ab = ba = 0.$$

As $\overline{Q}_1$ is generated over $\mathbb{F}$ by $1_1, 1_2, a, b$, modulo the relations $ab = 0 = ba$, it has dimension 4. (Here and throughout the paper we do not list obvious relations such as $a^2 = 0$, $1_1^2 = 1_1$, $a1_1 = a$, etc.)

**3.11.** For subsets $I$ and $J$ of $\mathbb{N}$ and an $I \cup J$-tuple $\underline{\tau} = (\tau_i)_{i \in I \cup J}$ of commuting automorphisms of $\mathbb{F}$, let $\mathcal{B}(\mathbb{F}, I, J, \underline{\tau})$ be the $\mathbb{F}$-category with a unique object, say $\omega$, whose endomorphism algebra is a unital associative algebra over $\mathbb{F}$ generated by $a_i, b_i, c_j, c_j^{-1}$, $i \in I$, $j \in J$, subject to the relations: $a_i b_i = b_i a_i = 0$, $a_i \lambda = \tau_i(\lambda) a_i$, $\lambda b_i = b_i \tau_i(\lambda)$ for $i \in I$ and $\lambda \in \mathbb{F}$, $c_j c_j^{-1} = 1 = 1_\omega$, $c_j \lambda = \tau_j(\lambda) c_j$ for $j \in J$ and $\lambda \in \mathbb{F}$, and $u_s v_t = v_t u_s$ for $u, v \in \{a, b, c\}$ and all $s, t \in I \cup J$.

**Proposition 3.12.** *Let $A_n$ be the nth Weyl algebra and assume $\mathcal{O}$ is an orbit in $D = \mathbb{K}[t_1, \ldots, t_n]$ under the automorphism group $\mathcal{G}$ (which is assumed to have a maximal break if $\mathcal{O}$ is degenerate and $n = \infty$). Let $\mathfrak{m}$ be the unique element of $\mathfrak{B}_\mathcal{O}$ in the nondegenerate case, and the designated maximal break of $\mathcal{O}$ in the degenerate case.*

1. *If $\operatorname{char} \mathbb{K} = 0$, then $\mathcal{S}_\mathcal{O} \cong \mathcal{A}(D/\mathfrak{m}, \underline{I}(\mathfrak{m}))$;*

2. *If $\operatorname{char} \mathbb{K} = p > 0$, then $\mathcal{S}_\mathcal{O} \cong \mathcal{B}(D/\mathfrak{m}, \underline{I}(\mathfrak{m}), \underline{I}(\mathfrak{m})^c, \underline{\tau}(\mathfrak{m}))$.*

**Proof.** Corollary 3.9 shows that the category $\mathcal{S}_\mathcal{O}$ with objects $\mathfrak{B}_\mathcal{O}$ is a skeleton of $\mathcal{C}_\mathcal{O}$.

Assume that $\operatorname{char} \mathbb{K} = 0$. When $\underline{I}(\mathfrak{m}) \neq \emptyset$, define the functor $G : \mathcal{A}(D/\mathfrak{m}, \underline{I}(\mathfrak{m})) \to \mathcal{S}_\mathcal{O}$ as follows:

$$G(\alpha) = \sigma^\alpha(\mathfrak{m}), \qquad G(a_{\alpha,i}) = X_{\sigma^\alpha(\mathfrak{m}),i}, \qquad G(b_{\alpha,i}) = Y_{\sigma^\alpha(\mathfrak{m}),i} \qquad (10)$$

where $\sigma^\alpha = \prod_i \sigma_i^{\alpha_i}$. From (3.3) it is easy to see that this is an isomorphism. Now suppose that $\underline{I}(\mathfrak{m}) = \emptyset$. The orbit is nondegenerate, and in this case, the functor $G : \mathcal{A}(D/\mathfrak{m}, \emptyset) \to \mathcal{S}_\mathcal{O}$ defined by

$$G(\omega) = \mathfrak{m}, \qquad G(1_\omega) = 1_\mathfrak{m} \qquad (11)$$

is an isomorphism.

We suppose that $\operatorname{char} \mathbb{K} = p$, and for $1 \leq i \leq n$, let $r_i$ denote the minimal positive integer such that $\sigma_i^{r_i}(\mathfrak{m}) = \mathfrak{m}$. Thus, $r_i = p$ or $1$. Then we define

$$\begin{aligned}
G : \mathcal{B}(D/\mathfrak{m}, \underline{I}(\mathfrak{m}), \underline{I}(\mathfrak{m})^c, \underline{\tau}(\mathfrak{m})) &\to \mathcal{S}_\mathcal{O} \\
G(\omega) &= \mathfrak{m} \\
G(a_i) &= X_{\sigma_i^{r_i-1}(\mathfrak{m}),i} X_{\sigma_i^{r_i-2}(\mathfrak{m}),i} \cdots X_{\sigma_i(\mathfrak{m}),i} X_{\mathfrak{m},i} \\
G(b_i) &= Y_{\mathfrak{m},i} Y_{\sigma_i(\mathfrak{m}),i} \cdots Y_{\sigma_i^{r_i-1}(\mathfrak{m}),i} \\
G(c_j) &= X_{\sigma_j^{r_j-1}(\mathfrak{m}),j} X_{\sigma_j^{r_j-2}(\mathfrak{m}),j} \cdots X_{\sigma_j(\mathfrak{m}),j} X_{\mathfrak{m},j}.
\end{aligned} \qquad (12)$$

It can be easily checked that this determines an isomorphism. □

**Remark 3.13.** The algebras $\mathcal{A}(D/\mathfrak{m}, \underline{I}(\mathfrak{m}))$ are finite-dimensional when $|\underline{I}(\mathfrak{m})| < \infty$; whereas the algebras $\mathcal{B}(D/\mathfrak{m}, \underline{I}(\mathfrak{m}), \underline{I}(\mathfrak{m})^c, \underline{\tau}(\mathfrak{m}))$ are always infinite-dimensional.

## 4 Simple weight modules for Weyl algebras.

The purpose of this section is to describe the simple weight modules for the Weyl algebra $A_n$ with $1 \leq n \leq \infty$. These results generalize those obtained in [BB], where the underlying field $\mathbb{K}$ was assumed to be algebraically closed. The case of the first Weyl algebra over an arbitrary field has been treated previously in for example ([DGO], [B1], [B2] and [BO1]).

**4.1.** Proposition 3.12 allows us to focus on the algebras of the form $\mathcal{A}(\mathbb{F}, I)$ in the characteristic 0 case. Thus, we assume $\mathbb{F}$ is a field of characteristic 0. For each $\alpha \in \mathfrak{ob}\mathcal{A}(\mathbb{F}, I)$, define a simple $\mathcal{A}(\mathbb{F}, I)$-module $S_\alpha$ such that $S_\alpha(\beta) = \delta_{\alpha,\beta}\mathbb{F}$ for all objects $\beta \in \mathfrak{ob}\mathcal{A}(\mathbb{F}, I)$, and let all morphisms be trivial.

In order to deal with the case $n = \infty$ we will need the following well-known result, which can be found in [DOF].

**Lemma 4.2.** *Let $\mathcal{C}$ be a category, and assume $M$ is a simple $\mathcal{C}$-module and $M(\alpha) \neq 0$ for some $\alpha \in \mathfrak{ob}\mathcal{C}$. Then $M(\alpha)$ is simple as a $\mathcal{C}(\alpha, \alpha)$-module. Conversely for any simple $\mathcal{C}(\alpha, \alpha)$-module $N$, there exists a unique (up to isomorphism) simple $\mathcal{C}$-module $M$ such that $M(\alpha) \cong N$ as $\mathcal{C}(\alpha, \alpha)$-modules.*

**Proposition 4.3.** *Any simple module over the algebra $\mathcal{A} = \mathcal{A}(\mathbb{F}, I)$ is isomorphic to $S_\alpha$ for some object $\alpha \in \mathfrak{ob}\mathcal{A}(\mathbb{F}, I)$.*

**Proof.** If $n < \infty$ and $I \subseteq \{1, \ldots, n\}$, then the statement is obvious since $\mathcal{A}$ is finite-dimensional. Suppose $n = \infty$, and let $M$ be a simple $\mathcal{A}$-module. For an object $\alpha \in \mathfrak{ob}\mathcal{A}$ denote by $M(\alpha)$, the module $M$ at $\alpha$. Assume that $M(\alpha) \neq 0$. By Lemma 4.2, $M(\alpha)$ is a simple $\mathcal{A}(\alpha, \alpha)$-module. Suppose $u \in \mathcal{A}(\alpha, \alpha)$, $u \neq 0, 1_\alpha$, and $u = (u_1)_{\beta^1, i_1} \cdots (u_m)_{\beta^m, i_m}$ where $u_j \in \{a, b\}$, $\beta^j \in \mathfrak{ob}\mathcal{A}$, and $i_j \in I$ for each $j = 1, \ldots, m$. Since $u \in \mathcal{A}(\alpha, \alpha)$, then for each $j$, there exists $k$ such that $\{u_j, u_k\} = \{a, b\}$ and $i_j = i_k$. Using the second set of relations (the commuting relations) in $\mathcal{A}(\mathbb{F}, I)$, we can rewrite $u$ as product $u = \prod_j a_{\beta^j, s_j} b_{\beta^j, s_j}$. But then the first set of relations in $\mathcal{A}(\mathbb{F}, I)$ implies that $u = 0$, which contradicts our assumption. Hence, $\mathcal{A}(\alpha, \alpha) = \mathbb{F}1_\alpha$, and $M(\alpha)$ is one-dimensional because it is simple as an $\mathcal{A}(\alpha, \alpha)$-module.

Consider the set $J$ of all $\beta \in \mathfrak{ob}\mathcal{A}$ such that $\alpha_j = \beta_j$ for all $j \in I$ except one and $M(\beta) \neq 0$. When $\beta \in J$ and $\beta_i \neq \alpha_i$, we say that $\beta$ is an $\alpha$-sink if $a_{\alpha,i} \neq 0$ (hence $b_{\alpha,i} = 0$) and is an $\alpha$-source if $b_{\alpha,i} \neq 0$ (hence $a_{\alpha,i} = 0$). Let $J_1$ (resp. $J_2$) be the set of all $\alpha$-sinks in $J$ (resp. the set of all $\alpha$-sources in $J$) and $J_3 = J \setminus (J_1 \cup J_2)$. Consider the following submodule of $M$: $M' = \mathcal{A}\big(\sum_{\beta \in J_1 \cup J_3} M(\beta)\big)$. We claim that $M' \cap M(\alpha) = 0$. Indeed, if $u \in \mathcal{A}$ and $uM(\beta) = M(\alpha)$ for some $\beta \in J_1 \cup J_3$, then it follows from the commuting relations that $u$ contains $b_{\alpha,i}$, where $\beta_i = 1$ and $\beta_j = \alpha_j$ for $j \neq i$. But because $\beta \in J_1 \cup J_3$, it must be that $b_{\alpha,i} = 0$, and hence $u = 0$. We conclude that $M' \cap M(\alpha) = 0$, so that $M' = 0$ by simplicity. This argument shows that whenever $M(\alpha) \neq 0$, then there are no $\alpha$-sinks. But then by simplicity, we must have $M = M(\alpha)$, and hence $M = S_\alpha$ in this case. □

**4.4.** In the characteristic $p > 0$ case, it suffices to consider the algebras $\mathcal{B}(\mathbb{F}, I, J, \underline{\tau})$ where $I, J \subseteq \mathbb{N}$ by Proposition 3.12. Assume $\Gamma \subseteq I$, and $\xi : \Gamma \to \{0, 1\}$. Associated to this data is the subalgebra

$$\mathfrak{R}_{\Gamma,\xi} = \mathbb{F}[d_i, c_j^{\pm 1}, \underline{\tau} \mid i \in \Gamma, j \in J] \tag{13}$$

of $\mathcal{B}(\mathbb{F}, I, J, \underline{\tau})$ where $d_i = a_i$ if $\xi(i) = 0$ and $d_i = b_i$ if $\xi(i) = 1$. (In (13), only the automorphisms $\tau_i$ with $i \in \Gamma \cup J$ are used, and $\mathfrak{R}_{\Gamma,\xi}$ is a skew polynomial algebra.) Then for a maximal ideal $\mathfrak{N}$ in $\mathfrak{R}_{\Gamma,\xi}$, define the simple module $S_{\Gamma,\xi,\mathfrak{N}} = \mathfrak{R}_{\Gamma,\xi}/\mathfrak{N}$. Note that when $I = \emptyset$, then $\Gamma = \emptyset$; and whenever $\Gamma = \emptyset$, we assume there is just one $\xi$ and $\mathfrak{R}_{\Gamma,\xi} := \mathbb{F}[c_j^{\pm 1}, \underline{\tau} \mid j \in J]$. Again in this case, for every maximal ideal $\mathfrak{N}$ of $\mathfrak{R}_{\Gamma,\xi}$, we define $S_{\Gamma,\xi,\mathfrak{N}} = \mathfrak{R}_{\Gamma,\xi}/\mathfrak{N}$.

**Proposition 4.5.** *The modules $S_{\Gamma,\xi,\mathfrak{N}}$ constitute an exhaustive list of pairwise nonisomorphic simple $\mathcal{B}(\mathbb{F}, I, J, \underline{\tau})$-modules.*

**Proof.** As the case $I = \emptyset$ is apparent, we may assume $S$ is a simple $\mathcal{B}(\mathbb{F}, I, J, \underline{\tau})$-module where $I \neq \emptyset$. It is clear that both $a_i S$ and $b_i S$ are submodules for any $i \in I$. Because $S$ is simple, each of them coincides either with $S$ or with $0$. But $a_i S = S$ implies $b_i S = 0$ and vice versa. Therefore, either $a_i S = b_i S = 0$, or $a_i S = S$, $b_i S = 0$, or $a_i S = 0$, $b_i S = S$ for any $i \in I$. Let $\Gamma = \{i \in I \mid \text{either } a_i S \neq 0 \text{ or } b_i S \neq 0\}$. Define $\xi : \Gamma \to \{0, 1\}$ by $\xi(i) = 0$ if $a_i S \neq 0$, $\xi(i) = 1$ if $b_i S \neq 0$. Then $\mathcal{B}(\mathbb{F}, I, J, \underline{\tau})$ modulo the annihilator of $S$ in $\mathbb{F}[a_i, b_i, \underline{\tau} \mid i \in I]$ is isomorphic to $\mathfrak{R}_{\Gamma, \xi}$. The result then follows. $\square$

**4.6.** Next we will describe the simple modules for the Weyl algebra $A = A_n$ corresponding to Propositions 4.3 and 4.5. As before, we assume $D = \mathbb{K}[t_1, \ldots, t_n]$ where $t_i = \partial_i x_i$.

Assume first that char $\mathbb{K} = 0$ and that $\mathcal{O}$ is a nondegenerate orbit of $\mathcal{G}$ on $\mathfrak{max}\, D$. Recall in this case that $\mathfrak{B}_{\mathcal{O}} = \{\mathfrak{m}\}$. Set

$$S(\mathcal{O}) = \bigoplus_{\mathfrak{n} \in \mathcal{O}} D/\mathfrak{n} \tag{14}$$

and define a left $A$-module structure on $S(\mathcal{O})$ by specifying for $i = 1, \ldots, n$

$$x_i(d + \mathfrak{n}) := \sigma_i(d) + \sigma_i(\mathfrak{n}), \qquad \partial_i(d + \mathfrak{n}) := t_i \sigma_i^{-1}(d) + \sigma_i^{-1}(\mathfrak{n}). \tag{15}$$

As $S(\mathcal{O})$ is generated by $1 + \mathfrak{m}$, we have that $S(\mathcal{O}) \cong A/A\mathfrak{m}$ where $1 + \mathfrak{m} \mapsto 1 + A\mathfrak{m}$.

Now suppose that char $\mathbb{K} = 0$, $\mathcal{O}$ is degenerate, and $\mathfrak{m}$ is the fixed maximal break. For $\mathfrak{p} \in \mathfrak{B}_{\mathcal{O}}$ set

$$S(\mathcal{O}, \mathfrak{p}) := \bigoplus_{\mathfrak{n} \in \mathcal{O}_{\mathfrak{p}}} D/\mathfrak{n}, \tag{16}$$

where $\mathcal{O}_{\mathfrak{p}}$ is as in (9). One can define a structure of a left $A$-module on $S(\mathcal{O}, \mathfrak{p})$ by the same formulae as in (15), but when the image is not in $S(\mathcal{O}, \mathfrak{p})$, the result is 0. Assuming $\mathfrak{p} = \sigma_{i_1}^{\delta_1} \cdots \sigma_{i_s}^{\delta_s}(\mathfrak{m})$, where $\underline{I}(\mathfrak{m}) = \{i_1, \ldots, i_s\}$, we have in this case $S(\mathcal{O}, \mathfrak{p}) \cong A/A(\mathfrak{p}, Z_{i_1}, \ldots, Z_{i_s})$ where $Z_k = x_k$ if $\mathfrak{p}$ is a break with respect to $k$, and $Z_k = \partial_k$ otherwise. The isomorphism is given by $1 + \mathfrak{p} \mapsto 1 + A(\mathfrak{p}, Z_{i_1}, \ldots, Z_{i_s})$. It follows from the construction that $S(\mathcal{O})$ and $S(\mathcal{O}, \mathfrak{p})$ are simple $A$-modules.

Suppose now that char $\mathbb{K} = p > 0$. Let $\mathfrak{m}$ be the designated maximal ideal when the orbit $\mathcal{O}$ is nondegenerate, and let $\mathfrak{m}$ be the fixed maximal break when $\mathcal{O}$ is degenerate. Assume $\underline{I} = \underline{I}(\mathfrak{m})$ is the set of breaks of $\mathfrak{m}$, and let $\Gamma$ be a subset of $\underline{I}$. Choose $\xi : \Gamma \to \{0, 1\}$. Then for $\mathfrak{n} \in \mathcal{O}$, define

$$\mathfrak{R}_{\Gamma,\xi}(\mathfrak{n}) = (D/\mathfrak{n})[d_i, c_j^{\pm 1}, \underline{\tau}(\mathfrak{n}) \mid i \in \Gamma, j \in \underline{I}^c], \tag{17}$$

where $d_i := a_i$ if $\xi(i) = 0$ and $d_i := b_i$ if $\xi(i) = 1$ as in (4.4). Note that if $\Gamma = \emptyset$, then $\mathfrak{R}_{\Gamma,\xi}(\mathfrak{n}) = (D/\mathfrak{n})[c_j^{\pm 1}, \underline{\tau}(\mathfrak{n}) \mid j = 1, \ldots, n]$ (and there is only one possible $\xi$). Let $\mathfrak{N}$ be a maximal ideal in the ring $\mathfrak{R}_{\Gamma,\xi}(\mathfrak{m})$.

Set

$$S(\mathcal{O}, \Gamma, \xi, \mathfrak{N}) := \bigoplus_{\mathfrak{n} \in \mathcal{O}} \mathfrak{R}_{\Gamma,\xi}(\mathfrak{n})/\sigma_{\mathfrak{n}}(\mathfrak{N}), \tag{18}$$

where $\sigma_{\mathfrak{n}}(\mathfrak{N})$ means apply $\sigma_{\mathfrak{n}}$ to the coefficients of elements of $\mathfrak{N}$. Define a structure of a left $A$-module on $S(\mathcal{O}, \Gamma, \xi, \mathfrak{N})$ by specifying for each $i \in \{1, 2, \ldots, n\}$,

$$x_i(f + \sigma_{\mathfrak{n}}(\mathfrak{N})) = \begin{cases} \sigma_i(f) + \sigma_i \sigma_{\mathfrak{n}}(\mathfrak{N}) & \text{if } \mathfrak{n} \neq \sigma \mathfrak{m} \text{ for any } \sigma \in \widehat{\mathcal{G}}_i; \\ b_i \sigma_i \sigma_{\mathfrak{n}}(f) + \sigma_i \sigma_{\mathfrak{n}}(\mathfrak{N}) & \text{if } \mathfrak{n} = \sigma \mathfrak{m} \text{ for some } \sigma \in \widehat{\mathcal{G}}_i, \\ & \quad i \in \Gamma, \text{ and } \xi(i) = 1; \\ c_i \sigma_i \sigma_{\mathfrak{n}}(f) + \sigma_i \sigma_{\mathfrak{n}}(\mathfrak{N}) & \text{if } \mathfrak{n} = \sigma \mathfrak{m} \text{ for some } \sigma \in \widehat{\mathcal{G}}_i \\ & \quad \text{and } i \in \underline{I}(\mathfrak{m})^c; \\ 0 & \text{in all other cases;} \end{cases} \tag{19}$$

$$\partial_i(f + \sigma_{\mathfrak{n}}(\mathfrak{N})) = \begin{cases} t_i \sigma_i^{-1} \sigma_{\mathfrak{n}}(f) + \sigma_i^{-1} \sigma_{\mathfrak{n}}(\mathfrak{N}) & \text{if } \mathfrak{n} \neq \sigma_i \sigma \mathfrak{m} \text{ for any } \sigma \in \widehat{\mathcal{G}}_i; \\ a_i t_i \sigma_i^{-1} \sigma_{\mathfrak{n}}(f) + \sigma_i^{-1} \sigma_{\mathfrak{n}}(\mathfrak{N}) & \text{if } \mathfrak{n} = \sigma_i \sigma \mathfrak{m} \text{ for some} \\ & \quad \sigma \in \widehat{\mathcal{G}}_i, \, i \in \Gamma, \& \, \xi(i) = 0; \\ c_i^{-1} t_i \sigma_i^{-1} \sigma_{\mathfrak{n}}(f) + \sigma_i^{-1} \sigma_{\mathfrak{n}}(\mathfrak{N}) & \text{if } \mathfrak{n} = \sigma_i \sigma \mathfrak{m} \text{ for some} \\ & \quad \sigma \in \widehat{\mathcal{G}}_i \text{ and } i \in \underline{I}(\mathfrak{m})^c, \\ 0 & \text{in all other cases;} \end{cases} \tag{20}$$

where $\widehat{\mathcal{G}}_i$ is the subgroup of $\mathcal{G}$ generated by the $\sigma_j$ for all $j \neq i$. It follows from the construction that $S(\mathcal{O}, \Gamma, \xi, \mathfrak{N})$ is simple $A$-module.

The next theorem provides a classification of simple weight $A_n$-modules, $1 \leq n \leq \infty$.

**Theorem 4.7.** *Let $A_n$ be the nth Weyl algebra for $1 \leq n \leq \infty$ over a field $\mathbb{K}$,*

(i) If $\operatorname{char} \mathbb{K} = 0$, then the modules $S(\mathcal{O})$ and $S(\mathcal{O}, \mathfrak{p})$, where $\mathfrak{p} \in \mathfrak{B}_{\mathcal{O}}$, constitute an exhaustive list of pairwise non-isomorphic simple weight $A_n$-modules with support in $\mathcal{O}$.

(ii) If $\operatorname{char} \mathbb{K} = p > 0$ then the modules $S(\mathcal{O}, \Gamma, \xi, \mathfrak{N})$ constitute an exhaustive list of pairwise non-isomorphic simple weight $A_n$-modules with support $\mathcal{O}$.

(In (i) and (ii), the orbit $\mathcal{O}$ is assumed to have a maximal break when $\mathcal{O}$ is degenerate and $n = \infty$.)

**Proof.** We fix an orbit $\mathcal{O}$ of $\mathbb{K}[t_1, \ldots, t_n]$ under the group $\mathcal{G}$ and let $\mathcal{C}_{\mathcal{O}}$ and $\mathcal{S}_{\mathcal{O}}$ be the associated category and skeleton as (3.3) and (3.5). We can identify $\mathcal{S}_{\mathcal{O}}$ with an algebra $\mathcal{A}(D/\mathfrak{m}, \underline{I}(\mathfrak{m}))$ or $\mathcal{B}(D/\mathfrak{m}, \underline{I}(\mathfrak{m}), \underline{I}(\mathfrak{m})^c, \underline{\tau}(\mathfrak{m}))$ using the functor $G$ defined by (10), (11), or (12) in the proof of Proposition 3.12. We will use the determination of modules for those algebras in Propositions 4.3 and 4.5. We denote by $E$ the equivalence functor $\mathcal{C}_{\mathcal{O}} \otimes_{\mathcal{S}_{\mathcal{O}}} -$ and let $F'$ be given by (8).

(i) $\operatorname{char} \mathbb{K} = 0$: When the orbit $\mathcal{O}$ is nondegenerate, it is easy to see that $F'EG(S_\omega) \cong S(\mathcal{O})$. When the orbit is degenerate, and $\mathfrak{m}$ is the designated maximal break with break set $\underline{I}(\mathfrak{m})$, we have for each $\alpha \in \mathcal{A}(D/\mathfrak{m}, \underline{I}(\mathfrak{m}))$ that $F'EG(S_\alpha) \cong S(\mathcal{O}, \sigma^\alpha(\mathfrak{m}))$ where $\sigma^\alpha = \prod_{i \in \underline{I}(\mathfrak{m})} \sigma_i^{\alpha_i}$.

(ii) $\operatorname{char} \mathbb{K} = p$: For the simple modules $S_{\Gamma, \xi, \mathfrak{N}}$ for $\mathcal{B}(D/\mathfrak{m}, \underline{I}(\mathfrak{m}), \underline{I}(\mathfrak{m})^c, \underline{\tau}(\mathfrak{m}))$ we have that $F'EG(S_{\Gamma, \xi, \mathfrak{N}}) \cong S(\mathcal{O}, \Gamma, \xi, \mathfrak{N})$. The statement in this case is a consequence of Proposition 4.5. $\square$

## 5 Examples.

**5.1. Case $A_1$**

Let $A_1 = \mathbb{K}[t](\sigma, t)$ be the first Weyl algebra, where $t = \partial x$ and $\sigma(t) = t - 1$. Simple weight $A_1$-modules have been described in [Bl], [B1], [B2] and [BO1], and they are the following:

1. Assume $\operatorname{char} \mathbb{K} = 0$, $\mathcal{O}$ is a nondegenerate orbit of $\mathfrak{max}(\mathbb{K}[t])$ under $\sigma$, and $\mathfrak{m} \in \mathcal{O}$. Then $\mathfrak{m}$ is a principal ideal generated by an irreducible polynomial different from $t$ over $\mathbb{K}$ and $S(\mathcal{O}) = A/A\mathfrak{m}$ is the corresponding simple $A$-module.

2. Assume $\operatorname{char} \mathbb{K} = 0$, $\mathcal{O}$ is degenerate, and $\mathfrak{m} \in \mathcal{O}$ is the fixed maximal break (which is unique here). Then $\mathfrak{m} = (t)$, $\mathfrak{B}_{\mathcal{O}} = \{\mathfrak{m}, \sigma(\mathfrak{m})\}$, $S(\mathcal{O}, \mathfrak{m}) \cong A/Ax$ and $S(\mathcal{O}, \sigma(\mathfrak{m})) \cong A/A\partial$.

3. Assume $\operatorname{char} \mathbb{K} = p > 0$, $\mathcal{O}$ is nondegenerate, and $\mathfrak{m}$ is the fixed maximal ideal in $\mathcal{O}$. If $\sigma(\mathfrak{m}) \neq \mathfrak{m}$, then $\mathfrak{R}_{\Gamma,\xi}(\mathfrak{m}) = (\mathbb{K}[t]/\mathfrak{m})[c^{\pm 1}]$ (here $\Gamma = \emptyset$ and there is just one $\xi$). Then the simple $A$-modules $S(\mathcal{O}, \emptyset, \xi, \mathfrak{N})$ are parametrized by the maximal ideals $\mathfrak{N} = (f)$ of $(\mathbb{K}[t]/\mathfrak{m})[c^{\pm 1}]$ generated by irreducible polynomials $f$ of $(\mathbb{K}[t]/\mathfrak{m})[c]$ different from $c$. When $\sigma(\mathfrak{m}) = \mathfrak{m}$, (for example, when $\mathfrak{m}$ is generated by a polynomial of the form $f(t) = t^p - t - \nu$ for some nonzero $\nu \in \mathbb{K}$ and $f$ is irreducible), then $\mathfrak{R} := \mathfrak{R}_{\Gamma,\xi}(\mathfrak{m}) = (\mathbb{K}[t]/\mathfrak{m})[c^{\pm 1}, \sigma]$ is a skew polynomial algebra. In this case, the simple $A$-modules $S(\mathcal{O}, \emptyset, \xi, \mathfrak{N})$ are parametrized by the maximal ideals $\mathfrak{N}$ of $\mathfrak{R}$. Any such maximal ideal is principal, $\mathfrak{N} = (f)$, where $f$ is irreducible in $\mathfrak{R}$, and $\mathfrak{R}/(f) \cong \mathfrak{R}/(g)$ if and only if $f$ and $g$ are similar. (See for example, [DGO].)

4. Assume $\operatorname{char} \mathbb{K} = p > 0$, $\mathcal{O}$ is degenerate, and $\mathfrak{m} \in \mathcal{O}$ is the fixed maximal break. Then $\mathfrak{m} = (t)$, and $\mathfrak{R}_{\Gamma,\xi}(\mathfrak{m}) = \mathbb{K}[t]/\mathfrak{m} \cong \mathbb{K}$ if $\Gamma = \emptyset$ and $\mathfrak{R}_{\Gamma,\xi}(\mathfrak{m}) \cong \mathbb{K}[d]$ if $\Gamma \neq \emptyset$. Hence in the first case, there is a unique simple module; while in the second, the simple modules $S(\mathcal{O}, \Gamma, \xi, \mathfrak{N})$ are parametrized by maximal ideals $\mathfrak{N} \subset \mathbb{K}[d]$.

It follows from (3) and (4) that if $\mathbb{K}$ is an algebraically closed field of characteristic $p$, then each irreducible module $S(\mathcal{O}, \Gamma, \xi, \mathfrak{N})$ has dimension $p$, as $\mathcal{O}$ has $p$ elements and each $\mathfrak{R}_{\Gamma,\xi}(\mathfrak{n})/\sigma_\mathfrak{n}(\mathfrak{N})$ is one-dimensional (see [C] for another approach to these modules). This need not be true if $\mathbb{K}$ is not algebraically closed. For example, if $\mathbb{K} = \mathbb{Z}_2$, and $\mathfrak{m}$ is the maximal ideal of $\mathbb{K}[t]$ generated by the irreducible polynomial $t^2 + t + 1$, then $\sigma(\mathfrak{m}) = \mathfrak{m}$ so $\mathcal{O} = \{\mathfrak{m}\}$. The maximal ideal $\mathfrak{N}$ of $(\mathbb{K}[t]/\mathfrak{m})[c^{\pm 1}, \sigma]$ generated by the irreducible polynomial $c^3 + c + 1$ gives an irreducible $A_1$-module $S(\mathcal{O}, \emptyset, \xi, \mathfrak{N})$ of dimension 6 over $\mathbb{Z}_2$.

### 5.2. Case $A_2$

Here we describe all simple weight $A$-modules for the second Weyl algebra $A = A_2$. If $\mathbb{K}$ is algebraically closed and $\operatorname{char} \mathbb{K} = 0$, then any simple weight $A$-module is a tensor product of two simple $A_1$-modules (cf. [BO2]). But in general this is not the case. We view $A$ as a generalized Weyl algebra $A = \mathbb{K}[t_1, t_2]((\sigma_1, \sigma_2), (t_1, t_2))$ as above, where $t_i = \partial_i x_i$, $\sigma_i(t_j) = \delta_{i,j}(t_i - 1)$ for $i, j = 1, 2$.

(i) Suppose that $\operatorname{char} \mathbb{K} = 0$, $\mathcal{O}$ is degenerate, and $\mathfrak{m} \in \mathcal{O}$ is the fixed maximal break. Then we have the following three possibilities.

1. $\mathfrak{m} = (t_1, t_2)$ and $\mathfrak{B}_\mathcal{O} = \{\mathfrak{m}, \sigma_1(\mathfrak{m}), \sigma_2(\mathfrak{m}), \sigma_1\sigma_2(\mathfrak{m})\}$,

$$S(\mathcal{O}, \mathfrak{m}) \cong A/A(x_1, x_2),$$
$$S(\mathcal{O}, \sigma_1(\mathfrak{m})) \cong A/A(\partial_1, x_2),$$
$$S(\mathcal{O}, \sigma_2(\mathfrak{m})) \cong A/A(x_1, \partial_2),$$
$$S(\mathcal{O}, \sigma_1\sigma_2(\mathfrak{m})) \cong A/A(\partial_1, \partial_2).$$

2. $\mathfrak{m} = (t_1, f)$, where $f \in \mathbb{K}[t_2]$ is an irreducible polynomial such that $t_2 \notin (f)$, and $\mathfrak{B}_\mathcal{O} = \{\mathfrak{m}, \sigma_1(\mathfrak{m})\}$,
$$S(\mathcal{O}, \mathfrak{m}) \cong A/A(f, x_1),$$
$$S(\mathcal{O}, \sigma_1(\mathfrak{m})) \cong A/A(f, \partial_1).$$

3. $\mathfrak{m} = (f, t_2)$, where $f \in \mathbb{K}[t_1]$ is an irreducible polynomial such that $t_1 \notin (f)$, and $\mathfrak{B}_\mathcal{O} = \{\mathfrak{m}, \sigma_2(\mathfrak{m})\}$,
$$S(\mathcal{O}, \mathfrak{m}) \cong A/A(f, x_2),$$
$$S(\mathcal{O}, \sigma_2(\mathfrak{m})) \cong A/A(f, \partial_2).$$

Note that in all cases the simple modules are tensor products of simple modules for $A_1$.

(ii) Assume now that $\mathcal{O}$ is a nondegenerate orbit in the characteristic $0$ case, and $\mathfrak{m} \subset \mathbb{K}[t_1, t_2]$ is the fixed maximal ideal in $\mathcal{O}$. Then $t_1 \notin \mathfrak{m}$ and $t_2 \notin \mathfrak{m}$, and $S(\mathcal{O}) = A/A\mathfrak{m}$ is the corresponding simple $A$-module.

(iii) Assume $\operatorname{char} \mathbb{K} = p > 0$, and $\mathfrak{m}$ is the designated maximal ideal in $\mathcal{O}$ when $\mathcal{O}$ is nondegenerate, and $\mathfrak{m}$ is the fixed maximal break with break set $\underline{I} = \underline{I}(\mathfrak{m})$ when $\mathcal{O}$ is degenerate. Then the simple $A$-modules $S(\mathcal{O}, \Gamma, \xi, \mathfrak{N})$ are parametrized by subsets $\Gamma \subseteq \underline{I}$, maps $\xi : \Gamma \to \{0, 1\}$, and maximal ideals $\mathfrak{N}$ of $\mathfrak{R}_{\Gamma, \xi}(\mathfrak{m}) = (\mathbb{K}[t_1, t_2]/\mathfrak{m})[d_i, c_j^{\pm 1}, \underline{\tau}(\mathfrak{m}) \mid i \in \Gamma, j \in \underline{I}^c]$ where $d_i = a_i$ if $\xi(i) = 0$ and $d_i = b_i$ if $\xi(i) = 1$. This leads to the following possibilities.

1. If $|\underline{I}| = 0$ and $\Gamma = \emptyset$, then the simple $A$-modules $S(\mathcal{O}, \emptyset, \xi, \mathfrak{N})$ are parametrized by maximal ideals $\mathfrak{N} \subset (\mathbb{K}[t_1, t_2]/\mathfrak{m})[c_j^{\pm 1}, \underline{\tau}(\mathfrak{m}) \mid j = 1, 2]$.

2. If $|\underline{I}| = 1$ and $\Gamma = \emptyset$, then $\mathfrak{R}_{\emptyset, \xi}(\mathfrak{m}) \cong (\mathbb{K}[t_1, t_2]/\mathfrak{m})[c_j^{\pm 1}, \underline{\tau}(\mathfrak{m})]$, and the modules $S(\mathcal{O}, \emptyset, \xi, \mathfrak{N})$ are parametrized by maximal ideals $\mathfrak{N} \subset (\mathbb{K}[t_1, t_2]/\mathfrak{m})[c_j^{\pm 1}, \underline{\tau}(\mathfrak{m})]$, where $\mathfrak{N} \neq (c_j)$.

3. If $|\underline{I}| = 1$ and $\Gamma = \underline{I} = \{i\}$, then $\mathfrak{m} = (t_i, f)$ for some polynomial $f = f(t_j) \notin \mathbb{K} t_j$ for $j \neq i$. The simple $A$-modules $S(\mathcal{O}, \Gamma, \xi, \mathfrak{N})$ are parametrized by maximal ideals $\mathfrak{N} \subset \left(\mathbb{K}[t_j]/(f)\right)[d_i][c_j^{\pm 1}, \underline{\tau}(\mathfrak{m})]$ and maps $\xi : \Gamma \to \{0, 1\}$ (which determine the $A$-action in (19) and (20)).

4. If $|\underline{I}| = 2$ and $\Gamma = \emptyset$, then $\mathfrak{m} = (t_1, t_2)$, $\mathfrak{R}_{\Gamma,\xi} \cong \mathbb{K}[t_1, t_2]/\mathfrak{m} \cong \mathbb{K}$, and the simple $A$-modules $S(\mathcal{O}, \emptyset, \xi, \mathfrak{N})$ are parametrized by maximal ideals $\mathfrak{N} \subset \mathbb{K}$, that is, by $\mathfrak{N} = 0$.

5. If $|\underline{I}| = 2$ and $\Gamma = \{i\}$, then $\mathfrak{R}_{\Gamma,\xi} \cong \mathbb{K}[d_i]$, and the corresponding simple $A$-modules $S(\mathcal{O}, \Gamma, \xi, \mathfrak{N})$ are parametrized by maximal ideals $\mathfrak{N} \subset \mathbb{K}[d_i]$ and maps $\xi$.

6. If $|\underline{I}| = 2$ and $\Gamma = \underline{I}$, then $\mathfrak{R}_{\Gamma,\xi} \cong \mathbb{K}[d_1, d_2]$, and the simple $A$-modules $S(\mathcal{O}, \Gamma, \xi, \mathfrak{N})$ are parametrized by maximal ideals $\mathfrak{N} \subset \mathbb{K}[d_1, d_2]$ and maps $\xi$.

### 5.3. Case of separable ideals

If $\mathfrak{m} = (\mathfrak{m}_1, \ldots, \mathfrak{m}_n)$, where $\mathfrak{m}_i$ is a maximal ideal of $\mathbb{K}[t_i]$ for each $i$, then $\mathbb{K}[t_1, \ldots, t_n]/\mathfrak{m} \cong \mathbb{K}[t_1]/\mathfrak{m}_1 \otimes \cdots \otimes \mathbb{K}[t_n]/\mathfrak{m}_n$. In particular, this is always the case when $\mathbb{K}$ is algebraically closed (as in [BB]). It implies that a skeleton $S_\mathcal{O}$ for $A_n$ is a tensor product of corresponding skeletons for $A_1$. Hence, in particular, any simple weight module is a tensor product of simple $A_1$-modules if $\mathbb{K}$ is algebraically closed or if $\mathbb{K}$ has characteristic 0.

## 6 Simple modules over the Heisenberg algebra.

Let $\mathbb{K} = \mathbb{C}$, the complex numbers. Here we study modules for the Heisenberg algebra which are $\mathbb{Z}$-graded with infinite-dimensional homogeneous components and have a nonzero central charge. Let $H = \mathbb{C}c \oplus \bigoplus_{i \in \mathbb{Z}\setminus\{0\}} \mathbb{C}e_i$ be an infinite-dimensional Heisenberg Lie algebra, where $[e_i, e_j] = \delta_{i,-j}c$, $[e_j, c] = 0$ for all $i \geq 1$ and all $j$. Let $\mathcal{H}_\mathbb{Z}$ denote the category of all $\mathbb{Z}$-graded $H$-modules $V$ such that $V = \bigoplus_{i \in \mathbb{Z}} V_i$ and $e_i V_j \subset V_{i+j}$. If $V \in \mathcal{H}_\mathbb{Z}$ is simple, then $c$ acts as a scalar, called the *(central) charge* of $V$. All simple modules with a zero charge have been classified in [Ch]. Simple modules for which the charge is nonzero but $0 < \dim_\mathbb{C} V_i < \infty$ for at least one $i$ have been described in [F]. There are relatively few known examples of simple modules $V$ with a nonzero charge and $\dim_\mathbb{C} V_i = \infty$ for all $i$. Here we will construct a series of such simple modules using the classification of simple weight $A_\infty$-modules. Without loss of generality we can assume that a simple module $V \in \mathcal{H}_\mathbb{Z}$ has central charge 1. Because the universal enveloping algebra of $H$ modulo the ideal generated by $c - 1$ is isomorphic to $A_\infty$, any simple module $V \in \mathcal{H}_\mathbb{Z}$ with central charge 1 becomes a module for $A = A_\infty$ by identifying $\partial_i = e_i$ and $x_i = e_{-i}$ for all $i > 0$. Any orbit $\mathcal{O}$ of $A$ is linear. Suppose $\mathcal{O}$ is nondegenerate and $\mathfrak{m}$ is the designated maximal ideal, and consider the simple $A$-module $S(\mathcal{O}) = A/A\mathfrak{m} = \oplus_{\mathfrak{p} \in \mathcal{O}} D/\mathfrak{p}$, where $D = \mathbb{C}[t_1, t_2, \dots]$ and $D/\mathfrak{p} \cong \mathbb{C}$ for any $\mathfrak{p}$. Then $S(\mathcal{O})$ is a $\mathbb{Z}$-graded simple $H$-module,

$$S(\mathcal{O}) = \sum_{i \in \mathbb{Z}} S^{(-i)}(\mathcal{O}),$$

where

$$S^{(-i)}(\mathcal{O}) = \sum_{\ell=0}^{\infty} \sum_{\substack{i_1, \dots, i_\ell \in \mathbb{Z} \\ \sum_{k=1}^{\ell} k i_k = i,\ i_\ell \neq 0}} D/\sigma_1^{i_1} \dots \sigma_\ell^{i_\ell}(\mathfrak{m}).$$

The $\ell = 0$ summand is $D/\mathfrak{m}$ if $i = 0$ and is 0 otherwise. Clearly, $S^{(-i)}(\mathcal{O})$ is infinite-dimensional for each $i$. Note that $S(\mathcal{O})$ remains simple as a $\mathbb{Z}$-graded module.

# 7 Representation type and indecomposable modules for $A_n$, $1 \leq n \leq \infty$.

In this section we classify the tame blocks in the category of locally-finite weight modules for the finite and infinite Weyl algebras $A_n$ over a field $\mathbb{K}$. We also determine the indecomposable modules in the tame blocks. These results generalize those in [BB] and [DGO].

We will use the following concepts which can be found in [D].

**Definition 7.1.** *Let $\mathbb{F}$ be an algebraically closed field, and let $\mathbb{F}\langle x, y \rangle$ be the free associative algebra with 1 over $\mathbb{F}$ on two generators $x, y$.*

- *An $\mathbb{F}$-category $\mathcal{C}$ is called* wild *if there exists a $\mathcal{C} - \mathbb{F}\langle x, y \rangle$- bimodule $M$, free of finite rank as a right $\mathbb{F}\langle x, y \rangle$-module, such that the functor $M \otimes_{\mathbb{F}\langle x,y \rangle} -$ preserves indecomposability and isomorphism classes.*

- *An $\mathbb{F}$-category $\mathcal{C}$ is called* tame *if, for each dimension vector $\underline{d}$, there exist a localization $R = \mathbb{F}[x]_f$ with respect to some $f \in \mathbb{F}[x]$ and a finite number of $\mathcal{C} - R$-bimodules $B_1, \cdots, B_n$ such that each $B_j$ is free of finite rank as a right $R$-module, and such that every indecomposable $X \in \mathcal{C}$-𝔣𝔡𝔪𝔬𝔡 with $\dim_{\mathbb{F}} X(\alpha) = \underline{d}_\alpha$ is isomorphic to $B_j \otimes_R S$ for some $j$ and some simple $R$-module $S$.*

In the case of an arbitrary field $\mathbb{F}$ we will say that an $\mathbb{F}$-category $\mathcal{C}$ is *wild* (resp., *tame*) if it is wild (resp., tame) over the algebraically closure $\overline{\mathbb{F}}$ of $\mathbb{F}$.

We will consider weight $A_n$-modules $V$ which are *locally-finite*, i.e. all $V_\mathfrak{m}$ are finite-dimensional vector spaces over $D/\mathfrak{m}$ where $D = \mathbb{K}[t_1, \ldots, t_n]$ and $t_i = \partial_i x_i$ as before. We let $W^{lf}(A_n)$ be the category of locally-finite weight $A_n$-modules. Also for each subset $T \subset \mathfrak{max}D$ we denote by $\mathcal{W}_T^{lf}(A_n)$ the full subcategory in $\mathcal{W}^{lf}(A_n)$ consisting of all modules $V$ with $\mathfrak{supp}(V) \subset T$. The indecomposable locally-finite weight modules over the Weyl algebra $A_1$ have been described in [DGO] (see also [BB] for the case of $A_n$ ($n < \infty$) over an algebraically closed field and [C] for induced module constructions for Weyl algebras $A_n$ ($n < \infty$) and quantum Weyl algebras). Because of the category isomorphisms of Sections 3 and 4, we can focus on the categories $\mathcal{A}(\mathbb{F}, I)$ and $\mathcal{B}(\mathbb{F}, I, J, \underline{\tau})$ for $\mathbb{F}$ a field (which eventually we specialize to $D/\mathfrak{m}$) and consider locally-finite indecomposable modules for them.

**7.2.** It follows from ([BB], Sec. 2.6) that $\mathcal{B}(\mathbb{F}, I, J, \underline{\tau})$ is wild if $|I| + |J| > 1$, and $\mathcal{A}(\mathbb{F}, I)$ is wild for $|I| > 2$. Therefore, it is enough to describe the indecomposables for the categories $\mathcal{A}(\mathbb{F}, I)$ with $|I| \leq 2$ and for $\mathcal{B}(\mathbb{F}, I, J, \underline{\tau})$ with $|I| + |J| = 1$, which corresponds to the $n = 1$ case.

In (4.1) we have constructed the simple modules for $\mathcal{A}(\mathbb{F}, I)$. The category $\mathcal{A}(\mathbb{F}, \emptyset)$ has a unique object $\omega$ and morphism set $\mathbb{F}1_\omega$. Thus it is a simple algebra, and any finite-dimensional indecomposable module for $\mathcal{A}(\mathbb{F}, \emptyset)$ is simple.

### 7.3. Indecomposables for $\mathcal{A}(\mathbb{F}, I)$, $|I| = 1$

Recall that the category $\mathcal{A}(\mathbb{F}, I)$, $|I| = 1$, is isomorphic to the category $\overline{\mathcal{Q}}_1 = \mathbb{F}\mathcal{Q}_1/\mathcal{R}$ corresponding to the following quiver and relations:

$$\mathcal{Q}_1: \quad \begin{array}{c} 1 \quad a \\ \circ \rightleftarrows \circ \\ b \quad 2 \end{array} \qquad ab = ba = 0.$$

We denote by $S_i$, $i = 1, 2$, the simple modules for $\overline{\mathcal{Q}}_1$ (compare (4.1)). Now let $M_a$ (resp., $M_b$) be the $\overline{\mathcal{Q}}_1$-module such that $M_a(1) = M_b(1) = \mathbb{F}e_1$, $M_a(2) = M_b(2) = \mathbb{F}e_2$, where the action is given by $ae_1 := e_2$ and $be_2 := 0$ (resp., $ae_1 := 0$ and $be_2 := e_1$). The following proposition is evident.

**Proposition 7.4.** *The modules $S_1, S_2, M_a, M_b$ constitute an exhaustive list of pairwise nonisomorphic finite-dimensional indecomposable $\overline{\mathcal{Q}}_1$-modules.*

## 7.5. Indecomposables for $\mathcal{A}(\mathbb{F}, I)$, $|I| = 2$

It is easy to see that category $\mathcal{A}(\mathbb{F}, I)$, $|I| = 2$, is isomorphic to the category $\overline{\mathcal{Q}}_2 := \mathbb{F}\mathcal{Q}_2/\mathcal{R}$ corresponding to the following quiver and relations:

$$\mathcal{Q}_2: \quad \begin{array}{c} 1 \circ \xrightleftharpoons[b_1]{a_1} \circ 2 \\ b_0 \updownarrow a_0 \quad a_2 \updownarrow b_2 \\ b_3 \\ 0 \circ \xrightleftharpoons[a_3]{} \circ 3 \end{array} \qquad \begin{array}{l} a_i b_i = b_i a_i = 0, \text{ for } i \in \{0, 1, 2, 3\} \\ a_i a_j = b_\ell b_m \text{ for } i, j, \ell, m \in \{0, 1, 2, 3\}, \\ \text{whenever this is possible} \end{array}$$

We regard the objects $0, 1, 2, 3$ as the elements of $\mathbb{Z}_4 = \mathbb{Z}/4\mathbb{Z}$ and denote the corresponding simple modules by $S_i$, $i \in \mathbb{Z}_4$.

**7.6.** Now suppose $M_i$, $i \in \mathbb{Z}_4$, is the $\overline{\mathcal{Q}}_2$-module such that $M_i(j) = \mathbb{F}e_j$ for $j \in \mathbb{Z}_4$, where for $M_i$ the action is given by $a_i e_i = e_{i+1}$, $a_{i+1} e_{i+1} = e_{i+2}$, $b_{i-1} e_i = e_{i-1}$, $b_{i-2} e_{i-1} = e_{i-2}$, and $u_j e_k = 0$ for all other instances of $u \in \{a, b\}$ and $j, k \in \mathbb{Z}_4$.

**7.7.** For each $n \in \mathbb{N}$, $n > 1$, and $j \in \mathbb{Z}_4$ define the $\overline{\mathcal{Q}}_2$-module $M_{n,j,0}$ (resp., $M_{n,j,1}$) as follows. Consider $n$ elements $e_1, \ldots, e_n$. For $\ell \in \mathbb{Z}_4$, an $\mathbb{F}$-basis of the vector space $M_{n,j,0}(\ell)$ (resp. $M_{n,j,1}(\ell)$) is the set of $e_k$ such that $j + k - 1 \equiv \ell \pmod{4}$. The action of $a_\ell$ and $b_{\ell-1}$ on $M_{n,j,0}(\ell)$ (resp. $M_{n,j,1}(\ell)$) is given by the rules:

$$a_\ell e_k = \begin{cases} e_{k+1} & \text{if } \ell \text{ is even (resp., odd)}, k < n, \& j + k - 1 \equiv \ell \pmod{4}; \\ 0 & \text{otherwise}. \end{cases}$$

$$b_{\ell-1} e_k = \begin{cases} e_{k-1} & \text{if } \ell \text{ is even (resp., odd)}, k > 1 \& j + k - 1 \equiv \ell \pmod{4}; \\ 0, & \text{otherwise}. \end{cases}$$

In all other cases, we have $u_j e_k = 0$ for $u \in \{a, b\}$ and $j, k \in \mathbb{Z}_4$.

**7.8.** We denote by $\text{Irr}_0 \mathbb{F}[x]$ the set of monic irreducible polynomials $f \neq x$ in $\mathbb{F}[x]$, and let $\text{Ind}_0 \mathbb{F}[x] = \{f^n \mid f \in \text{Irr}_0 \mathbb{F}[x] \text{ and } n \in \mathbb{N}\}$. For each $f(x) = x^e + \zeta_e x^{e-1} + \cdots + \zeta_1 \in \text{Ind}_0 \mathbb{F}[x]$, define the $\overline{\mathcal{Q}}_2$-module $M_{f,1}$ (resp., $M_{f,2}$) as follows: Set $M_{f,1}(i) := \mathbb{F}^e$ (resp. $M_{f,2}(i) = \mathbb{F}^e$) for $i \in \mathbb{Z}_4$, and define

$$\begin{aligned}
M_{f,1}(a_0) &= M_{f,1}(a_2) = M_{f,1}(b_1) = I_e \\
M_{f,1}(b_0) &= M_{f,1}(b_2) = M_{f,1}(a_1) = M_{f,1}(a_3) = 0, \\
M_{f,1}(b_3) &= \mathfrak{F}_f, \\
M_{f,2}(b_0) &= M_{f,2}(b_2) = M_{f,2}(a_1) = I_e, \\
M_{f,2}(a_0) &= M_{f,2}(a_2) = M_{f,2}(b_1) = M_{f,2}(b_3) = 0, \\
M_{f,2}(a_3) &= \mathfrak{F}_f,
\end{aligned}$$

where $\mathfrak{F}_f$ is the Frobenius (companion) matrix corresponding to the polynomial $f$:

$$\mathfrak{F}_f = \begin{pmatrix} 0 & 0 & \cdots & 0 & -\zeta_1 \\ 1 & 0 & \cdots & 0 & -\zeta_2 \\ 0 & 1 & \cdots & 0 & \vdots \\ \vdots & \vdots & \ddots & \vdots & -\zeta_{e-1} \\ 0 & 0 & \cdots & 1 & -\zeta_e \end{pmatrix}$$

As a consequence of Proposition 3.3.1 in [BB] we have

**Proposition 7.9.** *The modules $S_i$, $M_i$, $M_{n,i,0}$, $M_{n,i,1}$, $M_{f,1}$, $M_{f,2}$ where $f \in \operatorname{Ind}_0 \mathbb{F}[x]$, $i \in \mathbb{Z}_4$, and $n \in \mathbb{N}, n > 1$, constitute an exhaustive list of pairwise nonisomorphic finite-dimensional indecomposable $\overline{\mathcal{Q}}_2$-modules.*

**Theorem 7.10.** *Let $\mathcal{O}$ be an orbit of $\mathcal{G}$ on $\mathfrak{max}\, D$ for the Weyl algebra $A = A_n$ ($1 \leq n \leq \infty$) over $\mathbb{K}$ (which is assumed to have a maximal break when $\mathcal{O}$ is degenerate and $n = \infty$).*

(i) *If $\operatorname{char} \mathbb{K} = 0$, then the category $\mathcal{W}_{\mathcal{O}}^{lf}(A)$ is tame if and only if either $\mathcal{O}$ is nondegenerate or the order of the maximal break is less than or equal to 2;*

(ii) *If $\operatorname{char} \mathbb{K} = p > 0$, then the category $\mathcal{W}_{\mathcal{O}}^{lf}(A)$ is tame if and only if $n = 1$.*

**Proof.**

(i) (char $\mathbb{K} = 0$)

"$\Rightarrow$" This statement follows from (7.2).

"$\Leftarrow$" If $\mathcal{O}$ is nondegenerate with fixed maximal ideal $\mathfrak{m}$, then $\mathcal{W}_{\mathcal{O}}^{lf}(A) \cong \mathcal{A}(D/\mathfrak{m}, \emptyset)$-fdmod $\cong (D/\mathfrak{m})$-fdmod by Propositions 3.4 and 3.12, hence $\mathcal{W}_{\mathcal{O}}^{lf}(A)$ is tame. If $\mathcal{O}$ is degenerate with designated maximal break $\mathfrak{m}$ having order $|\underline{I}| = 1$, then $\mathcal{W}_{\mathcal{O}}^{lf}(A) \cong \mathcal{A}(D/\mathfrak{m}, \underline{I})$-fdmod by Propositions 3.4 and 3.12, hence $\mathcal{W}_{\mathcal{O}}^{lf}(A)$ is tame by Proposition 7.4. If $\mathcal{O}$ is degenerate with maximal break $\mathfrak{m}$ of order 2, then $\mathcal{W}_{\mathcal{O}}^{lf}(A) \cong \mathcal{A}(D/\mathfrak{m}, \underline{I})$-fdmod, $|\underline{I}| = 2$, by Propositions 3.4 and 3.12, hence $\mathcal{W}_{\mathcal{O}}^{lf}(A)$ is tame by Proposition 7.9.

(ii) (char $\mathbb{K} = p > 0$)

"$\Rightarrow$" If $\mathcal{B}(\mathbb{F}, I, J, \underline{\tau})$ is tame, then $|I| + |J| = 1$ by (7.2). This corresponds to the case $n = 1$ which has been treated in [DGO] and Sec. 4.1 of [BB].

"$\Leftarrow$" This is a consequence of Theorem 5.7 in [DGO]. $\square$

### 7.11. Indecomposables for the tame blocks

By Theorem 7.10, $A_1$ is the only Weyl algebra in characteristic $p$ having tame blocks, and indecomposable modules in these blocks have been determined in ([BB], Sec. 4.1). Thus, in what follows we concentrate on the characteristic 0 case and construct weight indecomposable $A_n$-modules for the tame blocks, using Propositions 3.4 and 3.12, and the description of $\mathcal{A}(\mathbb{F}, I)$-modules in case $|I| \leq 2$.

**Definition 7.12.** *Assume* char $\mathbb{K} = 0$, *and the orbit $\mathcal{O}$ has a maximal break $\mathfrak{m}$ of order $1$ with respect to $i$.*

- $M(\mathcal{O}, \mathfrak{p}) := S(\mathcal{O}, \mathfrak{p})$ *for* $\mathfrak{p} \in \mathfrak{B}_\mathcal{O}$. *(see (16) for the definition);*

- $M(\mathcal{O}, x_i) := \bigoplus_{\mathfrak{n} \in \mathcal{O}} D/\mathfrak{n}$, *where for each $j$:*

$$x_j(d + \mathfrak{n}) := \sigma_j(d) + \sigma_j(\mathfrak{n}) \quad \text{and}$$

$$\partial_j(d + \mathfrak{n}) := \begin{cases} 0 & \text{if } j = i,\ \sigma_i^{-1}(\mathfrak{n}) \in \mathcal{O}_\mathfrak{m}, \\ & \text{and } \mathfrak{n} \in \mathcal{O}_{\sigma_i(\mathfrak{m})}; \\ t_j \sigma_j^{-1}(d) + \sigma_j^{-1}(\mathfrak{n}) & \text{otherwise.} \end{cases}$$

- $M(\mathcal{O}, \partial_i) := \bigoplus_{\mathfrak{n} \in \mathcal{O}} D/\mathfrak{n}$, where for each $j = 1, \ldots, n$:

$$x_j(d + \mathfrak{n}) := \begin{cases} 0 & \text{if } j = i, \, \mathfrak{n} \in \mathcal{O}_{\mathfrak{m}}, \\ & \text{and } \sigma_i(\mathfrak{n}) \in \mathcal{O}_{\sigma_i(\mathfrak{m})}; \\ \sigma_j(d) + \sigma_j(\mathfrak{n}) & \text{otherwise,} \end{cases}$$

and

$$\partial_j(d + \mathfrak{n}) := t_j \sigma_i^{-1}(d) + \sigma_j^{-1}(\mathfrak{n}).$$

**Remark 7.13.** Let $\mathcal{O}, \mathfrak{m}, i$ be as in Definition 7.12. It is easy to see that

- $M(\mathcal{O}, \mathfrak{m}) \cong A/A(\mathfrak{m}, x_i)$;
- $M(\mathcal{O}, \sigma_i(\mathfrak{m})) \cong A/A(\sigma_i(\mathfrak{m}), \partial_i)$;
- $M(\mathcal{O}, x_i) \cong A/A\mathfrak{m}$;
- $M(\mathcal{O}, \partial_i) \cong A/A\sigma_i(\mathfrak{m})$.

**7.14.** Now suppose $\mathcal{O}$ is an orbit of $D = \mathbb{K}[t_1, \ldots, t_n]$ under the automorphism group $\mathcal{G}$ of $A_n$ having maximal break $\mathfrak{m}$ with respect to $i, j \in \{1, \ldots, n\}$ ($i \neq j$). Assuming $\mathbb{F} = D/\mathfrak{m}$, define a map $\gamma : \mathfrak{B}_{\mathcal{O}} \to \mathbb{Z}_4$ by the following rule: $\gamma(\mathfrak{m}) = 0$, $\gamma(\sigma_i(\mathfrak{m})) = 3$, $\gamma(\sigma_j(\mathfrak{m})) = 1$ and $\gamma(\sigma_i \sigma_j(\mathfrak{m})) = 2$.

To $M$ a $\overline{\mathcal{Q}}_2$-module, we associate a corresponding weight module $\mathfrak{M}$ for $A_n$ according to the following procedure: For $\mathfrak{p} \in \mathfrak{B}_{\mathcal{O}}$ suppose $\mathfrak{M}(\mathfrak{p}) \cong M(\gamma(\mathfrak{p}))$ (as $D/\mathfrak{m}$-vector spaces); via the map that takes $v \in M(\gamma(\mathfrak{p}))$ to $v_{\mathfrak{p}} \in \mathfrak{M}(\mathfrak{p})$. For each $\mathfrak{n} \in \mathcal{O}_{\mathfrak{p}}$ (compare (9)), suppose $\mathfrak{M}(\mathfrak{n}) \cong \mathfrak{M}(\mathfrak{p})$ via $v_{\mathfrak{p}} \mapsto v_{\mathfrak{n}}$ for all $v \in M(\gamma(\mathfrak{p}))$.

As $\sigma_{\mathfrak{n}}^{-1}$ induces an isomorphism from $D/\mathfrak{n}$ to $D/\mathfrak{m}$, we can consider $\mathfrak{M}(\mathfrak{n})$ as a $D/\mathfrak{n}$-module; hence as a $D$-module that $\mathfrak{n}$ annihilates. Therefore $\mathfrak{M} = \bigoplus_{\mathfrak{n} \in \mathcal{O}} \mathfrak{M}(\mathfrak{n})$ where $\mathfrak{M}(\mathfrak{n}) = \{u \in \mathfrak{M} \mid \mathfrak{n} u = 0\}$. For the $A_n$-action define:

- if $\mathfrak{n} \in \mathcal{O}_{\mathfrak{m}}$ and $\sigma_i(\mathfrak{n}) \in \mathcal{O}_{\sigma_i(\mathfrak{m})}$:

$$x_i v_{\mathfrak{n}} = (M(b_3)v)_{\sigma_i(\mathfrak{n})}, \qquad \partial_i v_{\sigma_i(\mathfrak{n})} = (M(a_3)v)_{\mathfrak{n}};$$
$$x_i v_{\sigma_j(\mathfrak{n})} = (M(a_1)v)_{\sigma_i \sigma_j(\mathfrak{n})}, \qquad \partial_i v_{\sigma_i \sigma_j(\mathfrak{n})} = (M(b_1)v)_{\sigma_j(\mathfrak{n})};$$

- if $\mathfrak{n} \in \mathcal{O}_{\mathfrak{m}}$ and $\sigma_j(\mathfrak{n}) \in \mathcal{O}_{\sigma_j(\mathfrak{m})}$:

$$x_j v_{\mathfrak{n}} = (M(a_0)v)_{\sigma_j(\mathfrak{n})}, \qquad \partial_j v_{\sigma_j(\mathfrak{n})} = (M(b_0)v)_{\mathfrak{n}};$$
$$x_j v_{\sigma_i(\mathfrak{n})} = (M(b_2)v)_{\sigma_i \sigma_j(\mathfrak{n})}, \qquad \partial_j v_{\sigma_i \sigma_j(\mathfrak{n})} = (M(a_2)v)_{\sigma_i(\mathfrak{n})};$$

- in all other cases:

$$x_\ell v_\mathfrak{n} = v_{\sigma_\ell(\mathfrak{n})}, \qquad \partial_\ell v_\mathfrak{n} = t_\ell v_{\sigma_\ell^{-1}(\mathfrak{n})}.$$

For $\mathfrak{p} \in \mathfrak{B}_\mathcal{O}$, we denote by $S(\mathcal{O}, \mathfrak{p})$, $M(\mathcal{O}, \mathfrak{p})$, $M(\mathcal{O}, n, \mathfrak{p}, 0)$, $M(\mathcal{O}, n, \mathfrak{p}, 1)$, $M(\mathcal{O}, f, 1)$, and $M(\mathcal{O}, f, 2)$ the weight $A_n$-module which corresponds via this process to the $\overline{\mathcal{Q}}_2$-module $S_{\gamma(\mathfrak{p})}$, $M_{\gamma(\mathfrak{p})}$, $M_{n,\gamma(\mathfrak{p}),0}$, $M_{n,\gamma(\mathfrak{p}),1}$, $M_{f,1}$ and $M_{f,2}$ respectively as in Proposition 7.9.

**Remark 7.15.** *Let $\mathcal{O}, \mathfrak{m}, i, j$, and $\mathfrak{p} \in \mathfrak{B}_\mathcal{O}$ be as in (7.14). It is easy to see that*

- $S(\mathcal{O}, \mathfrak{m}) \cong A/A(\mathfrak{m}, x_i, x_j)$;

- $S(\mathcal{O}, \sigma_i(\mathfrak{m})) \cong A/A(\sigma_i(\mathfrak{m}), \partial_i, x_j)$;

- $S(\mathcal{O}, \sigma_j(\mathfrak{m})) \cong A/A(\sigma_j(\mathfrak{m}), \partial_j, x_i)$;

- $S(\mathcal{O}, \sigma_i\sigma_j(\mathfrak{m})) \cong A/A(\sigma_i\sigma_j(\mathfrak{m}), \partial_i, \partial_j)$;

- $M(\mathcal{O}, \mathfrak{p}) \cong A/A\mathfrak{p}$.

**Theorem 7.16.** *Let $\mathcal{O}$ be an orbit of $\mathcal{G}$ on $\mathfrak{max}\, D$ for the Weyl algebra $A_n$ ($1 \leq n \leq \infty$) over a field $\mathbb{K}$ of characteristic 0 (which is assumed to have a maximal break if $\mathcal{O}$ is degenerate and $n = \infty$).*

- *If $\mathcal{O}$ is nondegenerate, then the simple module $S(\mathcal{O})$ (see Section 3 for definition), is the unique (up to isomorphism) indecomposable module in $\mathcal{W}_\mathcal{O}^{lf}(A_n)$;*

- *If $\mathcal{O}$ has a maximal break $\mathfrak{m}$ of order 1 with respect to $i$, then the modules $S(\mathcal{O}, \mathfrak{m})$, $S(\mathcal{O}, \sigma_i(\mathfrak{m}))$, $M(\mathcal{O}, x_i)$ and $M(\mathcal{O}, \partial_i)$, constitute an exhaustive list of pairwise nonisomorphic indecomposable modules in $\mathcal{W}_\mathcal{O}^{lf}(A_n)$;*

- *If $\mathcal{O}$ has a maximal break $\mathfrak{m}$ of order 2 with respect to $i$ and $j$, then the modules $S(\mathcal{O}, \mathfrak{p})$, $M(\mathcal{O}, \mathfrak{p})$, $M(\mathcal{O}, n, \mathfrak{p}, 0)$, $M(\mathcal{O}, n, \mathfrak{p}, 1)$, $M(\mathcal{O}, f, 1)$, and $M(\mathcal{O}, f, 2)$, where $f \in \mathrm{Ind}_0(D/\mathfrak{m})[x]$, $\mathfrak{p} \in \mathfrak{B}_\mathcal{O}$ and $n \in \mathbb{N}, n > 1$, constitute an exhaustive list of pairwise nonisomorphic indecomposable modules in $\mathcal{W}_\mathcal{O}^{lf}(A_n)$.*

**Proof.** As in the proof of Theorem 4.7, we use the functor $F'EG$, where $E$ is the equivalence functor $\mathcal{C}_\mathcal{O} \otimes_{\mathcal{S}_\mathcal{O}} -$, $F'$ is given by (8), and $G$ is defined by (10), (11), or (12).

(i) If the orbit $\mathcal{O}$ is nondegenerate, then it follows from Propositions 3.4 and 3.12, that $\mathcal{W}_\mathcal{O}^{lf}(A_n) \cong \mathcal{A}(D/\mathfrak{m}, \emptyset)\text{-fdmod} \cong D/\mathfrak{m}\text{-fdmod}$. Therefore every indecomposable module is simple in this case, and the statement follows from Theorem 7.10.

(ii) Assume $\mathcal{O}$ has a maximal break $\mathfrak{m}$ of order 1 with respect to $i$. Then it is easy to see that
$$F'EG(S_1) \cong S(\mathcal{O}, \mathfrak{m}), \qquad F'EG(S_2) \cong S(\mathcal{O}, \sigma_i(\mathfrak{m})),$$
$$F'EG(M_a) \cong M(\mathcal{O}, x_i), \qquad F'EG(M_b) \cong M(\mathcal{O}, \partial_i).$$
Consequently, the result in this case follows from Propositions 3.4, 3.12, and 7.4.

(iii) Suppose $\mathcal{O}$ has a maximal break $\mathfrak{m}$ of order 2 with respect to $i$ and $j$. Let $\gamma : \overline{\mathcal{Q}}_2 \to \mathfrak{B}_\mathcal{O}$ be as in (7.14). From
$$F'EG(S_k) \cong S(\mathcal{O}, \gamma^{-1}(k)), \qquad F'EG(M_k) \cong M(\mathcal{O}, \gamma^{-1}(k)),$$
$$F'EG(M_{n,k,\ell}) \cong M(\mathcal{O}, n, \gamma^{-1}(k), \ell), \qquad F'EG(M_{f,s}) \cong M_{\mathcal{O},f,s},$$
where $k \in \mathbb{Z}_4$, $\ell \in \{0,1\}$, and $s \in \{1,2\}$, we obtain the desired result from Propositions 3.4, 3.12, and 7.9. □

**Remark 7.17.** *Because the algebras $\mathcal{A}(D/\mathfrak{m}, \emptyset)$ and $\overline{\mathcal{Q}}_1$ are finite-dimensional and have only finitely many nonisomorphic finite-dimensional indecomposable modules, it follows from [A1], [A2] that any indecomposable module for those algebras is finite-dimensional. Therefore, in (i) and (ii) of Theorem 7.16 we obtain a classification of all indecomposable weight modules not just the locally-finite ones.*

## Acknowledgment

We gratefully acknowledge the support and hospitality of the following institutions:

- Universidade Federal do Rio Grande do Norte (the first author)

- the Mathematical Sciences Research Institute, Berkeley and the Fields Institute, Toronto (the second and third authors)

- the University of Wisconsin at Madison (the third author).

Departamento de Matemática, CCET
Universidade Federal do Rio Grande do Norte
Campus Universitário, Lagoa Nova
CEP 59072-970, Natal-RN, Brasil
e-mail address: bekkert@ccet.ufrn.br

Department of Mathematics
University of Wisconsin
Madison, WI 53706 USA
e-mail address: benkart@math.wisc.edu

Institute of Mathemática e Estatíśtica
Universidade de São Paulo
Caixa Postal 66281 São Paulo, CEP 05311-970 Brasil
e-mail address: futorny@ime.usp.br